\documentclass[11pt]{article}
\usepackage[left=1in,right=1in,top=1in,bottom=1in]{geometry}
\usepackage{times}
\usepackage{expl3}
\usepackage{cite}
\usepackage[table]{xcolor}
\usepackage{multirow}
\usepackage{stackengine} 
\usepackage{hhline}
\usepackage{lipsum}
\usepackage{titlesec}
\usepackage{wrapfig}
\usepackage{enumerate}
\usepackage{epsfig}
\usepackage{enumitem}
\usepackage{bbm}
\usepackage{calc}
\usepackage{graphicx}
\usepackage{amsmath}
\usepackage[title]{appendix}
\usepackage{amssymb}
\usepackage{epstopdf}
\usepackage{boldline}
\usepackage{calligra}
\usepackage{bm}
\usepackage{url}
\usepackage{blindtext}
\usepackage{accents}

\newcommand{\define}{\stackrel{\mbox{\tiny def}}{=}}

\newtheorem{definition}{Definition}
\newtheorem{theorem}{Theorem}
\newtheorem{corollary}{Corollary}
\newtheorem{lemma}{Lemma}

\newtheorem{example}{Example}
\usepackage{mathtools}
\usepackage{epstopdf}
\usepackage{balance}
\usepackage{thmtools}
\usepackage{thm-restate}
\usepackage{hyperref}
\usepackage{cleveref}

\usepackage[ruled,vlined]{algorithm2e}
\include{pythonlisting}
\newcommand{\ostar}{\mathbin{\mathpalette\make@circled\star}}

\makeatletter
\newcommand{\removelatexerror}{\let\@latex@error\@gobble}
\makeatother
\makeatletter
\newcommand*{\rom}[1]{\expandafter\@slowromancap\romannumeral #1@}
\makeatother

\ExplSyntaxOn
\newcommand\latinabbrev[1]{
  \peek_meaning:NTF . {
    #1\@}%
  { \peek_catcode:NTF a {
      #1.\@ }%
    {#1.\@}}}
\ExplSyntaxOff

\titleclass{\subsubsubsection}{straight}[\subsubsection]

\begin{document}
\vspace{1cm}
\title{Tail bounds for Multivariate Random Tensor Means}\vspace{1.8cm}
\author{Shih~Yu~Chang 
\thanks{Shih Yu Chang is with the Department of Applied Data Science,
San Jose State University, San Jose, CA, U. S. A. (e-mail: {\tt
shihyu.chang@sjsu.edu}).
           }}

\maketitle

\begin{abstract}
In our recent research endeavors, we have delved into the realm of tail bounds problems concerning bivariate random tensor means. In this context, tensors are treated as finite-dimensional operators. However, the longstanding challenge of extending the concept of operator means to scenarios involving more than two variables had persisted. The primary objective of this present study is to unveil a collection of tail bounds applicable to multivariate random tensor means. These encompass the weighted arithmetic mean, weighted harmonic mean, and the Karcher mean. These bounds are derived through the utilization of Ando-Hiai's inequalities, alongside tail bounds specifically tailored for multivariate random tensor means employing reverse Ando-Hiai's inequalities, which are rooted in Kantorovich constants. Notably, our methodology involves employing the concept of deformation for operator means with multiple variables, following the principles articulated in Hiai, Seo and Wada's recent work. Additionally, our research contributes to the expansion of the Karcher mean differentiable region from the vicinity of the diagonal identity element within the Cartesian product space of positive definite tensors to the vicinity of the general element within the Cartesian product space of positive definite tensors via the application of the inverse and implicit function theorem.
\end{abstract}


\section{Introduction}

Tensors are used in science and technologies for several reasons due to their versatility and ability to represent complex relationships and data structures. We discuss several important examples here. First, tensors are used to represent multidimensional data. Many real-world problems involve data that have multiple dimensions, such as time-series data, images, audio signals, and more. Tensors provide a natural and efficient way to represent and manipulate such data, making them essential in fields like image processing, computer vision, speech recognition, and signal processing. Second, tensors are fundamental building blocks in deep learning and neural networks. Neural networks process data in the form of tensors, and tensor-based operations, like matrix multiplications, are at the core of training and inference processes. This has led to remarkable advancements in various fields like computer vision, natural language processing, and robotics. Third, in physics and engineering, many systems are described by equations involving tensors. For example, in fluid dynamics, stress tensors describe the behavior of fluids under different conditions. Tensors are also used in electromagnetism, solid mechanics, and general relativity, among other disciplines, to represent physical quantities and their interactions. Fourth, in addition to deep learning, tensors are used in other machine learning techniques, like tensor decomposition methods and tensor-based data analysis. These methods are employed for data compression, dimensionality reduction, and feature extraction. Finally, tensors play a crucial role in robotics for representing spatial transformations and manipulator kinematics. They are used in control systems to model and control complex multi-input and multi-output systems. Overall, tensors offer a powerful and unified framework to handle multi-dimensional data and complex relationships, making them indispensable in various scientific and technological domains. Their ability to capture and process intricate information allows researchers from science and technology fields to solve challenging problems and make significant advancements in their respective fields~\cite{chang2023personalized,chang2023TLS,chang2022tPI,chang2022tPII,chang2022generaltail,chang2022TKF,
chang2022tensorq,chang2022randouble,chang2022TWF,chang2021convenient,chang2021tensorExp,chang2022generalized}.

In our recent works, we consider tail bounds problems for bivariate random tensor means, where tensors can be treated as finite dimensional operators~\cite{chang2023randomI,chang2023randomII}. However, the problem of extending the operator mean to more than two variables of matrices or operators had been a longstanding challenge.It was eventually resolved through two different approaches: an iteration method~\cite{ando2004geometric} and a Riemannian geometry method~\cite{bhatia2006riemannian}. Subsequently, the Riemannian geometry method has been significantly advanced by various authors, including~\cite{palfia2016operator}. In the present day, the multivariate geometric mean within the Riemannian geometric approach is commonly referred to as the \emph{Karcher mean} as it is determined by solving the so-called ``Karcher equation''. Additionally, another essential multivariate operator mean that has garnered recent active attention is the power mean, which was developed in~\cite{lawson2014karcher}. The purpose of this work is to present several tail bounds for multivariate random tensor means, e.g., weighted arithmetic mean, weighted harmonic mean, Karcher mean, based on Ando-Hiai’s inequalities and tail bounds for multivariate random tensor means with reverse Ando-Hiai’s inequalities derived from Kantorovich constants. The technique of deformation for $k$-variable operator means is used here according to~\cite{hiai2019ando, hiai2020ando}. Besides, we also extend Theorem 6.3 in~\cite{lawson2014karcher} about the Karcher mean differentiable region from the neighbor of the diagonal identity element in the Cartesian product space of positive definite tensors to the neighbor of general elements in the Cartesian product space of positive definite tensors by applying the inverse and implicit function theorem.

The rest of this paper is organized as follows. In Section~\ref{sec:Tensor Preliminaries}, we will review basic definitions about tensors used in this work, especially the notion about RTT. In Section~\ref{sec:Multivariate Random Tensor Means}, multivariate random tensor means and multivariate random tensor deformed means, and their associated examples like weighted arithmetic mean, weighted harmonic mean, Karcher mean, will be discussed. The continuity and differentiability of the Karcher mean will be discussed in Section~\ref{sec:Continuity and Differentiability of Karcher Means}. In Section~\ref{sec:Tail bounds for Multivariate Random Tensor Power Mean and Karcher Mean}, we will derive tail bounds for multivariate random tensor means. In Section~\ref{sec:Tail bounds for Multivariate Random Tensor Deformed Means Based on Ando-Hiai Type Inequalities}, we will generalize tail bounds from multivariate random tensor mean to multivariate random tensor deformed means. Finally, in Section~\ref{sec:Tail bounds for Multivariate Random Tensor Means Based on Reverse Ando-Hiai Type Inequalities}, we consider the reverse of Ando-Hiai Type inequalities via Kantorovich constants for multivariate random tensor means and multivariate random tensor deformed means. 

\noindent \textbf{Nomenclature:} The sets of complex and real numbers are denoted by $\mathbb{C}$ and $\mathbb{R}$, respectively. The set of natural numbers is represented by $\mathbb{N}$. A scalar is denoted by an either italicized or Greek alphabet such as $x$ or $\beta$; a vector is denoted by a lowercase bold-faced alphabet such as $\bm{x}$; a matrix is denoted by an uppercase bold-faced alphabet such as $\bm{X}$; a tensor is denoted by a calligraphic alphabet such as $\mathcal{x}$ or $\mathcal{X}$.

\section{Tensor Preliminaries}\label{sec:Tensor Preliminaries}

In this section, we will review required definitions and basic facts about tensors. The default product between two tensors $\mathcal{A} \in \mathbb{C}^{I_1 \times \dots \times I_N \times I_1 \times \dots \times I_N}$ and $\mathcal{B} \in \mathbb{C}^{I_1 \times \dots \times I_N \times I_1 \times \dots \times I_N}$ is Einstein product with order $N$, denoted by $\star_N$. We will specfiy the exact product symbol if the product is not $\star_N$.

As the matrix eigen-decomposition theorem is crucial in various linear algebra theory and applications, we will have a parallel decomposition theorem for Hermitian tensors. From Theorem 5.2 in~\cite{ni2019hermitian}, every Hermitian tensor $\mathcal{H} \in  \mathbb{C}^{I_1 \times \dots \times I_N \times I_1 \times \dots \times I_N}$ has the following decomposition:
\begin{eqnarray}\label{eq:Hermitian Eigen Decom}
\mathcal{H} &=& \sum\limits_{i=1}^r \lambda_i \mathcal{U}_i  \star_1 \mathcal{U}^{H}_i, \mbox{
~with~~$\langle \mathcal{U}_i, \mathcal{U}_i \rangle =1$ and $\langle \mathcal{U}_i, \mathcal{U}_j \rangle = 0$ for $i \neq j$,}
\end{eqnarray}
where $\lambda_i \in \mathbb{R}$ and $\mathcal{U}_i \in \mathbb{C}^{I_1 \times \dots \times I_N \times 1}$. Here tensors $\mathcal{U}_i$ are orthogonal tensors each other since $\langle \mathcal{U}_i, \mathcal{U}_j \rangle = 0$ for $i \neq j$. The values $\lambda_i$ are named as \emph{Hermitian eigenvalues}, and the minimum integer of $r$ to decompose a Hermitian tensor as in Eq.~\eqref{eq:Hermitian Eigen Decom} is called \emph{Hermitian tensor rank}. In this work, we assume that all Hermitian tensors discussed in this work are full rank, i.e., $r=\prod\limits_{j=1}^N I_j$. A \emph{positive definite} (PD) tensor is a Hermitian tensor with all \emph{Hermitian eigenvalues} are positive. A \emph{semipositive definite} (SPD) tensor is a Hermitian tensor of which all \emph{Hermitian eigenvalues} are nonnegative. The symbols $\lambda_{\min}(\mathcal{X})$ and $\lambda_{\max}(\mathcal{X})$ will represent the minimum and maximum eigenvalues of the tensor $\mathcal{X}$. Given two Hermitian tensors $\mathcal{X}$ and $\mathcal{Y}$, we write $\mathcal{X} \preceq \mathcal{Y}$ if $\mathcal{Y}-\mathcal{X}$ is a SPD tensor, and write $\mathcal{X} \succeq \mathcal{Y}$ if $\mathcal{X}-\mathcal{Y}$ is a SPD tensor. We adopt $S_{\mbox{\tiny PD}}$ to represent a set of PD tensors, and $S_{\mbox{\tiny PSD}}$ to represent a set of SPD tensors.

Let us represent the Hermitian eigenvalues of a Hermitian tensor $\mathcal{H} \in \mathbb{C}^{I_1 \times \dots \times I_N \times I_1 \times \dots \times I_N} $ in decreasing order by the vector $\vec{\lambda}(\mathcal{H}) = (\lambda_1(\mathcal{H}), \cdots, \lambda_r(\mathcal{H}))$, where $r$ is the Hermitian rank of the tensor $\mathcal{H}$. We use $\mathbb{R}_{\geq 0} (\mathbb{R}_{> 0})$ to represent a set of nonnegative (positive) real numbers. Let $\left\Vert \cdot \right\Vert_{\rho}$ be a unitarily invariant tensor norm, i.e., $\left\Vert \mathcal{H}\star_N \mathcal{U}\right\Vert_{\rho} = \left\Vert \mathcal{U}\star_N \mathcal{H}\right\Vert_{\rho} = \left\Vert \mathcal{H}\right\Vert_{\rho} $,  where $\mathcal{U}$ is for unitary tensor. Let $\rho : \mathbb{R}_{\geq 0}^r \rightarrow \mathbb{R}_{\geq 0}$ be the corresponding gauge function that satisfies H${\rm \ddot{o}}$lder’s inequality so that 
\begin{eqnarray}\label{eq:def gauge func and general unitarily invariant norm}
\left\Vert \mathcal{H} \right\Vert_{\rho} = \left\Vert |\mathcal{H}| \right\Vert_{\rho} = \rho(\vec{\lambda}( | \mathcal{H} | ) ),
\end{eqnarray}
where $ |\mathcal{H}|  \define \sqrt{\mathcal{H}^H \star_N \mathcal{H}} $. 

The notion about  \emph{random tensor topology} (RTT) is provided by the following definition.
\begin{definition}\label{def:conv in mean}
We say that a sequence of random tensor $\mathcal{X}_n \in \mathbb{C}^{I_1 \times \dots \times I_N \times I_1 \times \dots \times I_N}$ converges to the random tensor $\mathcal{X} \in \mathbb{C}^{I_1 \times \dots \times I_N \times I_1 \times \dots \times I_N}$ with respect to the tensor norm $\left\Vert \cdot \right\Vert_{\rho}$ in the sense of RTT, if we have
\begin{eqnarray}
\mathbb{E}\left( \left\Vert \mathcal{X}_n \right\Vert_{\rho} \right)~~~\mbox{exists,}
\end{eqnarray}
and
\begin{eqnarray}
\lim\limits_{n \rightarrow \infty}\mathbb{E}\left( \left\Vert \mathcal{X}_n  - \mathcal{X} \right\Vert_{\rho} \right) = 0. 
\end{eqnarray}
We adopt the notation $\lim\limits_{n \rightarrow \infty}\mathcal{X}_n = \mathcal{X}$ to represent that random tensors $\mathcal{X}_n$ converges to the random tensor $\mathcal{X}$ with respect to the tensor norm $\left\Vert \cdot \right\Vert_{\rho}$ in the sense of RTT.
\end{definition}

\section{Multivariate Random Tensor Means}\label{sec:Multivariate Random Tensor Means}

\subsection{Multivariate Random Tensor Means and Deformed Means}\label{sec:Multivariate Random Tensor Means and Deformed Means}

We have to introduce the Thompson metric first. The Thompson metric between two PD tensors $\mathcal{X}, \mathcal{Y} \in \mathbb{C}^{I_1 \times \dots \times I_N \times I_1 \times \dots \times I_N}$, represented by $d_{T}\left(\mathcal{X}, \mathcal{Y}\right) $ on PD tensors is defined by
\begin{eqnarray}\label{eq:Thompson metric def}
d_{T}\left(\mathcal{X}, \mathcal{Y}\right) \define \lambda_{\max}\left(\log\left(\mathcal{X}^{-1/2}\star_N \mathcal{Y}\star_N \mathcal{X}^{-1/2}\right)\right)=\log \max\{\bm{\alpha}(\mathcal{X}/\mathcal{Y}),\bm{\alpha}(\mathcal{Y}/\mathcal{X})\},
\end{eqnarray}
where $\bm{\alpha}(\mathcal{X}/\mathcal{Y}) \define \inf\{\alpha>0: \mathcal{X} \preceq \alpha \mathcal{Y}\}$~\cite{thompson1963certain}.

The concept of bivariate operator means was formulated as an axiomatic way by Kubo-Ando~\cite{kubo1980means}. Let $ \sigma: S_{\mbox{\tiny PD}}\times S_{\mbox{\tiny PD}} \rightarrow S_{\mbox{\tiny PD}}$ be a map of \emph{bivaraite tensor mean}, and all random tensors used in the following list with dimensions $\mathbb{C}^{I_1 \times \dots \times I_N \times I_1 \times \dots \times I_N}$, we extend this notion to random tensors as follows:
\begin{enumerate}[label=(\roman*)]\label{list:Kubo-Ando axioms 2 Arg}
\item Monotonicity: $\mathcal{X} \preceq \mathcal{Z}$, $\mathcal{Y}\preceq \mathcal{W}$ $\Longrightarrow$ 
$\mathcal{X}\sigma\mathcal{Y} \preceq \mathcal{Z}\sigma\mathcal{W}$, \mbox{~almost surely};
\item Transformer Inequality: 
\begin{eqnarray}
\mathcal{Z}\star_N (\mathcal{X}\sigma\mathcal{Y})\star_N \mathcal{Z}\preceq (\mathcal{Z}\star_N \mathcal{X}\star_N \mathcal{Z}) \sigma (\mathcal{Z}\star_N \mathcal{Y}\star_N \mathcal{Z}), \mbox{~almost surely},
\end{eqnarray}
where $\mathcal{Z} \in  S_{\mbox{\tiny SPD}}$;
\item Monotone continuity: $\mathcal{X}_n \searrow \mathcal{X}$,  $\mathcal{Y}_n \searrow \mathcal{Y}$,   $\Longrightarrow$ $\mathcal{X}_n \sigma \mathcal{Y}_n \searrow \mathcal{X} \sigma \mathcal{Y}$, where $\mathcal{X}_n \searrow \mathcal{X}$ indicates that $\mathcal{X}_1 \succeq \mathcal{X}_2 \succeq \cdots $ and $\lim\limits_{n \rightarrow \infty}\mathcal{X}_n = \mathcal{X}$ in RTT, on the other hand, $\mathcal{X}_n \nearrow \mathcal{X}$,  $\mathcal{Y}_n \nearrow \mathcal{Y}$,   $\Longrightarrow$ $\mathcal{X}_n \sigma \mathcal{Y}_n \nearrow \mathcal{X} \sigma \mathcal{Y}$, where $\mathcal{X}_n \nearrow \mathcal{X}$ indicates that $\mathcal{X}_1 \preceq \mathcal{X}_2 \preceq \cdots $ and $\lim\limits_{n \rightarrow \infty}\mathcal{X}_n = \mathcal{X}$ in RTT;
\item Normalized condition: $\mathcal{I} \sigma \mathcal{I} = \mathcal{I}$.
\end{enumerate}
Note that transformer inequality will become equal if the tensor $\mathcal{Z}$ belongs to $S_{\mbox{\tiny PD}}$ (invertible).

In~\cite{kubo1980means}, they showed that there is a one-to-one correspondence between the operator mean $\sigma$ and the non-negative operator monotone function $g_{\sigma}$ on $(0, \infty)$ with $g(1) = 1$ determined by 
\begin{eqnarray}
g(x)\mathcal{I} &=& \mathcal{I}\sigma(x \mathcal{I}),
\end{eqnarray}
where $x >0 $; and
\begin{eqnarray}\label{eq:sigma f relation}
\mathcal{X}\sigma \mathcal{Y} &=& \mathcal{X}^{-1/2}\star_N g(\mathcal{X}^{1/2}\star_N \mathcal{Y} \star_N \mathcal{X}^{-1/2})\star_N \mathcal{X}^{1/2},
\end{eqnarray}
where $\mathcal{X}, \mathcal{Y} \in \mathbb{C}^{I_1 \times \dots \times I_N \times I_1 \times \dots \times I_N}$ are PD tensors. Such corresponding function $g_{\sigma}$ is named as \emph{corresponding function} of $\sigma$. We say that the operator $\sigma$ is \emph{power monotone increasing} (p.m.i.) if $g_{\sigma}(x^p) \geq g^p_{\sigma}(x)$ for all $x > 0$ and $p \geq 1$. Using the limitation method given by~\cite{chang2023randomII}, we can extend Eq.~\eqref{eq:sigma f relation} from PD tensors to SPD tensors in RTT.

The bivariate Kubo-Ando axioms provided in Section~\eqref{list:Kubo-Ando axioms 2 Arg} can be extended to the multivariate version given in Section~\eqref{list:Kubo-Ando axioms mul Arg}.  Let $\mathfrak{M}: S^k_{\mbox{\tiny PD}} \rightarrow S_{\mbox{\tiny PD}}$ be a map of \emph{$k$-variable tensor mean}, and all random tensors used in the following list with dimensions $\mathbb{C}^{I_1 \times \dots \times I_N \times I_1 \times \dots \times I_N}$, the multivariate Kubo-Ando axioms is given by the follows:
\begin{enumerate}[label=(\Roman*)]\label{list:Kubo-Ando axioms mul Arg}
\item Monotonicity: If $\mathcal{X}_i, \mathcal{Y}_i \in S_{\mbox{\tiny PD}}$ and $\mathcal{X}_i \preceq \mathcal{Y}_i$ for $1 \leq i \leq k$, then we have $\mathfrak{M}(\mathcal{X}_1,\cdots,\mathcal{X}_k) \preceq \mathfrak{M}(\mathcal{Y}_1,\cdots,\mathcal{Y}_k)$, almost surely;
\item Congruence invariance: For every $\mathcal{X}_i \in S_{\mbox{\tiny PD}}$, where $1 \leq i \leq k$ and any invertible tensor $\mathcal{A}$, we have 
\begin{eqnarray}\label{eq:congruence inv}
\mathcal{A}^{\mathrm{H}}\star_N \mathfrak{M}(\mathcal{X}_1, \cdots, \mathcal{X}_k)\star_N \mathcal{A} 
&=& \mathfrak{M}(\mathcal{A}^{\mathrm{H}}\star_N \mathcal{X}_1 \star_N \mathcal{H}, \cdots, 
\mathcal{A}^{\mathrm{H}}\star_N \mathcal{X}_k \star_N \mathcal{H}),
\end{eqnarray}
where $\mathrm{H}$ is a Hermitan transpose;
\item Monotone continuity: $\mathcal{X}_i,  \mathcal{X}_{i,n} \in S^k_{\mbox{\tiny PD}}$ for $1 \leq i \leq k$ and $n \in \mathbb{N}$. If either $\mathcal{X}_{i,n} \searrow \mathcal{X}_i$, or $\mathcal{X}_{i,n} \nearrow \mathcal{X}_i$ as $n \rightarrow \infty$ for each $i$, we have
\begin{eqnarray}
\lim\limits_{n \rightarrow \infty}\mathfrak{M}(\mathcal{X}_{1,n}, \cdots, \mathcal{X}_{k,n})&=&
\mathfrak{M}(\mathcal{X}_{1}, \cdots, \mathcal{X}_{k}),
\end{eqnarray}
where the limitation is in the sense of RTT;
\item Normalized condition: $\mathfrak{M}(\mathcal{I},\cdots,\mathcal{I}) = \mathcal{I}$.
\end{enumerate}

Given $k$ PD tensors, $\mathcal{A}_1, \cdots, \mathcal{A}_k$, we consider the following tensor-valued equation
\begin{eqnarray}\label{eq:2.1}
\mathcal{X}&=&\mathfrak{M}(\mathcal{X}\sigma\mathcal{A}_1,\cdots,\mathcal{X}\sigma\mathcal{A}_k).
\end{eqnarray}
The solution $\mathcal{X}$ of Eq.~\eqref{eq:2.1} shall play an important role in determining the multivariate tensor mean with respect $\mathcal{A}_1, \cdots, \mathcal{A}_k$.  
 
Before presenting the main result of this section, we need following lemmas about the Thompson metric.
\begin{lemma}\label{lma:2.2}
For every $\mathcal{X}_i, \mathcal{Y}_i \in S_{\mbox{\tiny PD}}$, where $1 \leq i \leq k$, we have
\begin{eqnarray}\label{eq1:lma:2.2}
d_{T}(\mathfrak{M}(\mathcal{X}_1,\cdots,\mathcal{X}_k),\mathfrak{M}(\mathcal{Y}_1,\cdots,\mathcal{Y}_k))
\leq \max\limits_{1 \leq i \leq k}d_{T}(\mathfrak{M}(\mathcal{X}_i,\mathcal{Y}_i),~\mbox{almost surely.}
\end{eqnarray} 
\end{lemma}
\textbf{Proof:}
Let $\kappa>0$ be the largest term for $d_{T}(\mathfrak{M}(\mathcal{X}_i,\mathcal{Y}_i)$ among all $i$, we have
\begin{eqnarray}\label{eq2:lma:2.2}
e^{-\kappa}\mathcal{X}_i \preceq \mathcal{Y}_i \preceq e^{\kappa}\mathcal{X}_i,~\mbox{almost surely.}
\end{eqnarray} 
Then, from axioms (I) and (II) in the multivariate Kubo-Ando axioms given at Section~\eqref{list:Kubo-Ando axioms mul Arg}, we have 
\begin{eqnarray}\label{eq3:lma:2.2}
\mathfrak{M}(\mathcal{X}_1,\cdots,\mathcal{X}_k)&=&\mathfrak{M}(e^{-\kappa}\mathcal{X}_1,\cdots,e^{-\kappa}\mathcal{X}_k)\nonumber \\
&\preceq&\mathfrak{M}(\mathcal{Y}_1,\cdots,\mathcal{Y}_k),~\mbox{almost surely;}
\end{eqnarray} 
and
\begin{eqnarray}\label{eq4:lma:2.2}
\mathfrak{M}(\mathcal{Y}_1,\cdots,\mathcal{Y}_k)&\preceq&\mathfrak{M}(e^{\kappa}\mathcal{X}_1,\cdots,e^{\kappa}\mathcal{X}_k)\nonumber \\
&=&e^{\kappa}\mathfrak{M}(\mathcal{X}_1,\cdots,\mathcal{X}_k),~\mbox{almost surely.}
\end{eqnarray} 
Then, this Lemma is proved from Eqs.~\eqref{eq3:lma:2.2} and~\eqref{eq4:lma:2.2}.
$\hfill \Box$

The next Lemma shows that the Thompson metric will be reduced by taking the mean with respect to a PD tensor. 
\begin{lemma}\label{lma:2.3}
Given $\mathcal{X}, \mathcal{Y}, \mathcal{A} \in S_{\mbox{\tiny PD}}$ and $\mathcal{X} \neq \mathcal{Y}$, then 
$d_{T}(\mathcal{X}\sigma\mathcal{A},\mathcal{Y}\sigma\mathcal{A}) < d_{T}(\mathcal{X},\mathcal{Y})$.
\end{lemma}
\textbf{Proof:}
We wish to show that if $\mathcal{A} \prec \mathcal{B}$ is almost surely, we have $\mathcal{X}\sigma\mathcal{A} \prec \mathcal{X}\sigma\mathcal{B}$ is almost surely.  Let $g_{\sigma}$ be the one-to-one correspondence function with respect to the tensor mean $\sigma$, then we have 
\begin{eqnarray}\label{eq1:lma:2.3}
\mathcal{X}\sigma\mathcal{A}&=&\mathcal{X}^{1/2}g_{\sigma}(\mathcal{X}^{-1/2}\mathcal{A}\mathcal{X}^{-1/2})\mathcal{X}^{1/2}\nonumber \\
&\prec& \mathcal{X}^{1/2}g_{\sigma}(\mathcal{X}^{-1/2}\mathcal{B}\mathcal{X}^{-1/2})\mathcal{X}^{1/2}\nonumber \\
&=& \mathcal{X}\sigma\mathcal{B},
\end{eqnarray}
where we use that the function $g_{\sigma}$ is strictly tensor increasing on $(0, \infty)$. 

Set $\kappa \define d_{T}(\mathcal{X},\mathcal{Y}) >0$ and the fact that $e^{-\kappa}\mathcal{A} \prec \mathcal{A} \prec e^{\kappa}\mathcal{A}$ for any tensor $\mathcal{A}$, we have
\begin{eqnarray}\label{eq2:lma:2.3}
\mathcal{Y}\sigma\mathcal{A}\preceq e^{\kappa}\mathcal{X}\sigma\mathcal{A} \prec 
e^{\kappa}\mathcal{X}\sigma e^{\kappa}\mathcal{A}=e^{\kappa}(\mathcal{X}\sigma\mathcal{A}),
\end{eqnarray}  
and 
\begin{eqnarray}\label{eq3:lma:2.3}
\mathcal{Y}\sigma\mathcal{A}\succeq e^{-\kappa}\mathcal{X}\sigma\mathcal{A} \succ 
e^{-\kappa}\mathcal{X}\sigma e^{-\kappa}\mathcal{A}=e^{-\kappa}(\mathcal{X}\sigma\mathcal{A}).
\end{eqnarray}  
From Eqs.~\eqref{eq2:lma:2.3},~\eqref{eq3:lma:2.3} and Thompson meetric definition provided by Eq.~\eqref{eq:Thompson metric def}, this Lemma is proved.
$\hfill \Box$

We are ready to present the main result of this section.
\begin{theorem}\label{thm:2.1}
\begin{enumerate}[label=(\Roman*)]
\item Given $k$ PD tensors, $\mathcal{A}_1, \cdots, \mathcal{A}_k$, there exists a unique solution $\mathcal{X} \in S_{\mbox{\tiny PD}}$ satisfying the following:
\begin{eqnarray}\label{eq1:thm:2.1}
\mathcal{X}&=&\mathfrak{M}(\mathcal{X}\sigma\mathcal{A}_1,\cdots,\mathcal{X}\sigma\mathcal{A}_k);
\end{eqnarray}
\item Let $\mathfrak{M}_{\sigma}(\mathcal{A}_1,\cdots,\mathcal{A}_k)$ be the unique solution of Eq.~\eqref{eq1:thm:2.1}, then $\mathfrak{M}_{\sigma}: S^k_{\mbox{\tiny PD}} \rightarrow S_{\mbox{\tiny PD}}$ is a $k$-variable mean satisfying Kubo-Ando axioms given at Section~\eqref{list:Kubo-Ando axioms mul Arg};
\item If $\mathcal{X} \in S_{\mbox{\tiny PD}}$ and $\mathcal{X} \preceq \mathfrak{M}(\mathcal{X}\sigma\mathcal{A}_1,\cdots,\mathcal{X}\sigma\mathcal{A}_k)$, then $\mathcal{X} \preceq \mathfrak{M}_{\sigma}(\mathcal{A}_1,\cdots,\mathcal{A}_k)$;
\item If $\mathcal{X} \in S_{\mbox{\tiny PD}}$ and $\mathcal{X} \succeq \mathfrak{M}(\mathcal{X}\sigma\mathcal{A}_1,\cdots,\mathcal{X}\sigma\mathcal{A}_k)$, then $\mathcal{X} \succeq \mathfrak{M}_{\sigma}(\mathcal{A}_1,\cdots,\mathcal{A}_k)$.
\end{enumerate}
\end{theorem}
\textbf{Proof:}
For item (I) proof, we first select a constant $\alpha>0$ such that $\mathcal{A}_1, \cdots, \mathcal{A}_k \in S_{\alpha}$, where $S_{\alpha}$ is a set defined by $S_{\alpha} \define \{\mathcal{X} \in S_{\mbox{\tiny PD}}:\alpha^{-1}\mathcal{I} \preceq \mathcal{X} \preceq  \alpha\mathcal{I}\}$. We also define the map $G:S_{\mbox{\tiny PD}} \rightarrow S_{\mbox{\tiny PD}}$ by 
\begin{eqnarray}\label{eq2:thm:2.1}
G(\mathcal{X}) \define \mathfrak{M}_{\sigma}(\mathcal{A}_1,\cdots,\mathcal{A}_k).
\end{eqnarray}
From the multivariate Kubo-Ando axioms given at Section~\eqref{list:Kubo-Ando axioms mul Arg}, we have $G(\mathcal{X}) \in S_{\alpha}$ if $\mathcal{X} \in S_{\alpha}$ and monotonicity of $G(\mathcal{X})$. If we set $\mathcal{Y}_0 \define \alpha \mathcal{I}$, we have $\mathcal{Y}_0 \succeq G(\mathcal{Y}_0) \succeq G^2(\mathcal{Y}_0) \succeq \cdots \succeq \alpha^{-1}\mathcal{I}$, from (III) at the multivariate Kubo-Ando axioms given at Section~\eqref{list:Kubo-Ando axioms mul Arg}, we have $G^n(\mathcal{Y}_0) \searrow \mathcal{X}_0$ for some $\mathcal{X}_0 \in S_{\mbox{\tiny PD}}$. Then, we have
\begin{eqnarray}\label{eq3:thm:2.1}
\mathcal{X}_0 &=& \lim\limits_{n \rightarrow \infty}G(G^n(\mathcal{Y}_0)) \nonumber \\
&=& \lim\limits_{n \rightarrow \infty}\mathfrak{M}(G^n(\mathcal{Y}_0)\sigma\mathcal{A}_1,\cdots,G^n(\mathcal{Y}_0)\sigma\mathcal{A}_k)\nonumber \\
&=& \mathfrak{M}(\mathcal{X}_0\sigma\mathcal{A}_1,\cdots,\mathcal{X}_0\sigma\mathcal{A}_k),
\end{eqnarray}
where the limitation is taken in the sense of RTT. Therefore, $\mathcal{X}_0$ is the solution of Eq.~\eqref{eq1:thm:2.1}. 

For the uniqueness of the solution $\mathcal{X}_0$, if both solutions $\mathcal{X}_1$ and $\mathcal{X}_2$ with $\mathcal{X}_1 \neq \mathcal{X}_2$ satisfy the Eq.~\eqref{eq1:thm:2.1}, we have
\begin{eqnarray}\label{eq4:thm:2.1}
d_{T}(\mathcal{X}_1, \mathcal{X}_2) \leq \max\limits_{1 \leq i \leq k}d_{T}(\mathcal{X}_1\sigma\mathcal{A}_i,\mathcal{X}_2\sigma\mathcal{A}_i)<d_{T}(\mathcal{X}_1, \mathcal{X}_2), 
\end{eqnarray}
where Lemmas~\ref{lma:2.2} and~\ref{lma:2.3} are applied. Eq.~\eqref{eq4:thm:2.1} shows contradiction, and we should have a unique solution of Eq.~\eqref{eq1:thm:2.1}.

We shall prove items (III) and (IV) first. For item (III), since we assume that $\mathcal{X} \in S_{\mbox{\tiny PD}}$ and $\mathcal{X} \preceq \mathfrak{M}(\mathcal{X}\sigma\mathcal{A}_1,\cdots,\mathcal{X}\sigma\mathcal{A}_k)$, we have $\mathcal{X} \preceq G(\mathcal{X}) \preceq \cdots$. Again, by selecting proper $\alpha>0$, we have $\mathcal{X}, \mathcal{A}_1, \cdots, \mathcal{A}_n \in S_{\alpha}$. Then, we have $G^{n}(\mathcal{X}) \preceq \alpha \mathcal{I}$ for all $n$, which implies $G^n(\mathcal{X})\nearrow \mathcal{X}_0$ for some $\mathcal{X}_0 \in S_{\mbox{\tiny PD}}$ and $G(\mathcal{X}_0)=\mathcal{X}_0$. From the item (III) in the multivariate Kubo-Ando axioms given at Section~\eqref{list:Kubo-Ando axioms mul Arg}, we have $\mathcal{X} \preceq \mathcal{X}_0=\mathfrak{M}_{\sigma}(\mathcal{A}_1, \cdots, \mathcal{A}_k)$ as desired. 

For item (IV), since we assume that $\mathcal{X} \in S_{\mbox{\tiny PD}}$ and $\mathcal{X} \succeq \mathfrak{M}(\mathcal{X}\sigma\mathcal{A}_1,\cdots,\mathcal{X}\sigma\mathcal{A}_k)$, we have $\mathcal{X} \succeq G(\mathcal{X}) \succeq \cdots$. By selecting proper $\alpha>0$, we have $\mathcal{X}, \mathcal{A}_1, \cdots, \mathcal{A}_n \in S_{\alpha}$. Then, we have $G^{n}(\mathcal{X}) \succeq \alpha^{-1} \mathcal{I}$ for all $n$, which implies $G^n(\mathcal{X})\searrow \mathcal{X}_0$ for some $\mathcal{X}_0 \in S_{\mbox{\tiny PD}}$ and $G(\mathcal{X}_0)=\mathcal{X}_0$. From the item (III) in the multivariate Kubo-Ando axioms given at Section~\eqref{list:Kubo-Ando axioms mul Arg}, we have $\mathcal{X} \succeq \mathcal{X}_0=\mathfrak{M}_{\sigma}(\mathcal{A}_1, \cdots, \mathcal{A}_k)$ as we wish.

For item (II), suppose we are given two sets of random tensors $\mathcal{A}_i \in S_{\mbox{\tiny PD}}$ and $\mathcal{B}_i \in S_{\mbox{\tiny PD}}$ such that $\mathcal{A}_i \preceq \mathcal{B}_i$ for $1 \leq i \leq k$ and set $\mathcal{Y}_0 \define \mathfrak{M}_{\sigma}(\mathcal{B}_1,\cdots,\mathcal{B}_k)$. Then we have
\begin{eqnarray}\label{eq5:thm:2.1}
\mathcal{Y}_0&=&\mathfrak{M}(\mathcal{Y}_0\sigma\mathcal{B}_1,\cdots,\mathcal{Y}_0\sigma\mathcal{B}_k)\nonumber \\
&\succeq&\mathfrak{M}(\mathcal{Y}_0\sigma\mathcal{A}_1,\cdots,\mathcal{Y}_0\sigma\mathcal{A}_k)\nonumber \\
&\Longrightarrow_1& \\
\mathcal{Y}_0&\succeq&\mathfrak{M}(\mathcal{A}_1,\cdots,\mathcal{A}_k),
\end{eqnarray}
where we use the item (IV) in this Thoerem (proved above) in $\Longrightarrow_1$. Hence, we prove the item (I) at the multivariate Kubo-Ando axioms given at Section~\eqref{list:Kubo-Ando axioms mul Arg}. 

Let $\mathcal{X}_0 \define \mathfrak{M}_{\sigma}(\mathcal{A}_1,\cdots,\mathcal{A}_k)$, then from the item (II) at the multivariate Kubo-Ando axioms given at Section~\eqref{list:Kubo-Ando axioms mul Arg}, we have 
\begin{eqnarray}\label{eq6:thm:2.1}
\mathcal{Z}^{\mathrm{H}}\mathcal{X}_0\mathcal{Z}&=&
\mathcal{Z}^{\mathrm{H}}\mathfrak{M}(\mathcal{X}_0\sigma\mathcal{A}_1,\cdots,\mathcal{X}_0\sigma\mathcal{A}_k)\mathcal{Z} \nonumber \\
&=& \mathfrak{M}((\mathcal{Z}^{\mathrm{H}}\mathcal{X}_0\mathcal{Z})\sigma(\mathcal{Z}^{\mathrm{H}}\mathcal{A}_1\mathcal{Z}),\cdots,(\mathcal{Z}^{\mathrm{H}}\mathcal{X}_0\mathcal{Z})\sigma(\mathcal{Z}^{\mathrm{H}}\mathcal{A}_k\mathcal{Z})) \nonumber \\
&=& \mathfrak{M}_{\sigma}(\mathcal{Z}^{\mathrm{H}}\mathcal{A}_1\mathcal{Z},\cdots,\mathcal{Z}^{\mathrm{H}}\mathcal{A}_k\mathcal{Z}),
\end{eqnarray}
where $\mathcal{Z}$ is any invertible tensor. Therefore, we prove the item (II) at the multivariate Kubo-Ando axioms given at Section~\eqref{list:Kubo-Ando axioms mul Arg}. 

Let $\mathcal{A}_i, \mathcal{A}_{i,n} \in S_{\mbox{\tiny PD}}$ and assume that $\mathcal{A}_{i,n}  \searrow \mathcal{A}_{i}$ as $n \rightarrow \infty$ for $1 \leq i \leq k$. By setting, $\mathcal{X}_n \define \mathfrak{M}_{\sigma}(\mathcal{A}_{1,n},\cdots,\mathcal{A}_{k,n})$, from the fact (I) in this theorem (proved above), we have $\mathcal{X}_1 \succeq \mathcal{X}_2 \succeq \cdots$ and $\mathcal{X}_n \succeq \mathfrak{M}_{\sigma}(\mathcal{A}_1,\cdots,\mathcal{A}_k)$. Therefore, $\mathcal{X}_n \searrow \mathcal{X}_0$ for some $\mathcal{X}_0 \in S_{\mbox{\tiny PD}}$. Because $\mathcal{X}_n \sigma \mathcal{A}_{i,n} \searrow \mathcal{X}_0 \sigma \mathcal{A}_{i}$ for $1 \leq i \leq k$, from the (III) at the multivariate Kubo-Ando axioms given at Section~\eqref{list:Kubo-Ando axioms mul Arg} with respect to $\mathfrak{M}$, we have $\mathcal{X}_0 = \mathfrak{M}(\mathcal{X}_0\sigma\mathcal{A}_1,\cdots,\mathcal{X}_0\sigma\mathcal{A}_k)=\mathfrak{M}_{\sigma}(\mathcal{A}_1,\cdots,\mathcal{A}_k)$. This shows that $\mathfrak{M}_{\sigma}$ is downward continuous. On the other hand, let $\mathcal{B}_i, \mathcal{B}_{i,n} \in S_{\mbox{\tiny PD}}$ and assume that $\mathcal{B}_{i,n}  \nearrow \mathcal{A}_{i}$ as $n \rightarrow \infty$ for $1 \leq i \leq k$. By setting, $\mathcal{Y}_n \define \mathfrak{M}_{\sigma}(\mathcal{B}_{1,n},\cdots,\mathcal{B}_{k,n})$, from the fact (I) in this theorem (proved above), we have $\mathcal{Y}_1 \preceq \mathcal{Y}_2 \preceq \cdots$ and $\mathcal{Y}_n \preceq \mathfrak{M}_{\sigma}(\mathcal{B}_1,\cdots,\mathcal{B}_k)$. Therefore, $\mathcal{Y}_n \nearrow \mathcal{Y}_0$ for some $\mathcal{Y}_0 \in S_{\mbox{\tiny PD}}$. Because $\mathcal{Y}_n \sigma \mathcal{B}_{i,n} \searrow \mathcal{Y}_0 \sigma \mathcal{Y}_{i}$ for $1 \leq i \leq k$, from the (III) at the multivariate Kubo-Ando axioms given at Section~\eqref{list:Kubo-Ando axioms mul Arg} with respect to $\mathfrak{M}$, we have $\mathcal{Y}_0 = \mathfrak{M}(\mathcal{Y}_0\sigma\mathcal{B}_1,\cdots,\mathcal{Y}_0\sigma\mathcal{B}_k)=\mathfrak{M}_{\sigma}(\mathcal{B}_1,\cdots,\mathcal{B}_k)$. This shows that $\mathfrak{M}_{\sigma}$ is upward continuous. Therefore, we prove the item (III) at the multivariate Kubo-Ando axioms given at Section~\eqref{list:Kubo-Ando axioms mul Arg} with respect to $\mathfrak{M}_{\sigma}$.

Finally, due to that $\mathfrak{M}(\mathcal{I}\sigma\mathcal{I},\cdots,\mathcal{I}\sigma\mathcal{I})=\mathfrak{M}_{\sigma}(\mathcal{I},\cdots,\mathcal{I})=\mathcal{I}$, we prove the item (IV) at the multivariate Kubo-Ando axioms given at Section~\eqref{list:Kubo-Ando axioms mul Arg} with respect to $\mathfrak{M}_{\sigma}$.
$\hfill \Box$

Note that $\mathfrak{M}_{\sigma}$ in item (II) of Theorem~\ref{thm:2.1} is called the deformed mean of $\mathfrak{M}$ by $\sigma$ .

\subsection{Examples of Multivariate Means}\label{sec:Examples of Multivariate Means}

In this section, we will present several examples about multivariate means that satisfy the multivariate Kubo-Ando axioms given at Section~\eqref{list:Kubo-Ando axioms mul Arg}. 

\begin{example}\label{exp:adjoint}
Given a multivariate mean $\mathfrak{M}$, the \emph{adjoin}t of $\mathfrak{M}$, denoted by $\mathfrak{M}^*$ is defined by
\begin{eqnarray}\label{eq:adjoint def}
\mathfrak{M}^*&=&\mathfrak{M}^{-1}(\mathcal{A}_1^{-1},\cdots,\mathcal{A}_k^{-1}),
\end{eqnarray}
where $\mathcal{A}_i \in S_{\mbox{\tiny PD}}$ for $1 \leq i \leq k$. The multivariate mean $\mathfrak{M}^*$ satisfies the multivariate Kubo-Ando axioms given at Section~\eqref{list:Kubo-Ando axioms mul Arg}. 
\end{example}

\begin{example}\label{exp:weighted arithmetric mean}
Given a probability vector $\bm{w}=[w_1,\cdots,w_k]$ with length $k$ and $k$ random tensors $\mathcal{A}_i \in S_{\mbox{\tiny PD}}$ for $1 \leq i \leq k$, the weighted arithmetic mean with respect to $\bm{w}$, denoted as $\mathfrak{M}_{A,\bm{w}}(\mathcal{A}_1,\cdots,\mathcal{A}_k)$, can be defined by
\begin{eqnarray}
\mathfrak{M}_{A,\bm{w}}(\mathcal{A}_1,\cdots,\mathcal{A}_k) \define \sum\limits_{i=1}^k w_i \mathcal{A}_i.
\end{eqnarray}
The multivariate mean $\mathfrak{M}_{A,\bm{w}}(\mathcal{A}_1,\cdots,\mathcal{A}_k)$ satisfies the multivariate Kubo-Ando axioms given at Section~\eqref{list:Kubo-Ando axioms mul Arg}. 
\end{example}

\begin{example}\label{exp:weighted harmonic mean}
Given a probability vector $\bm{w}=[w_1,\cdots,w_k]$ with length $k$ and $k$ random tensors $\mathcal{A}_i \in S_{\mbox{\tiny PD}}$ for $1 \leq i \leq k$, the weighted harmonic mean with respect to $\bm{w}$, denoted as $\mathfrak{M}_{H,\bm{w}}(\mathcal{A}_1,\cdots,\mathcal{A}_k)$, can be defined by
\begin{eqnarray}
\mathfrak{M}_{H,\bm{w}}(\mathcal{A}_1,\cdots,\mathcal{A}_k) \define \left(\sum\limits_{i=1}^k w_i \mathcal{A}^{-1}_i\right)^{-1}.
\end{eqnarray}
The multivariate mean $\mathfrak{M}_{H,\bm{w}}(\mathcal{A}_1,\cdots,\mathcal{A}_k)$ satisfies the multivariate Kubo-Ando axioms given at Section~\eqref{list:Kubo-Ando axioms mul Arg}. 
\end{example}

The bivariate power tensor mean for two PD tensors $\mathcal{X} \in \mathbb{C}^{I_1 \times \dots \times I_N \times I_1 \times \dots \times I_N}$ and $\mathcal{Y} \in \mathbb{C}^{I_1 \times \dots \times I_N \times I_1 \times \dots \times I_N}$, represented by $\mathcal{X}\#_{q}\mathcal{Y}$, is defined by
\begin{eqnarray}\label{eq:bivaraite power mean def}
\mathcal{X}\#_{q}\mathcal{Y} \define \mathcal{X}^{1/2} \star_N (\mathcal{X}^{-1/2}\star_N\mathcal{Y}\star_N\mathcal{X}^{-1/2})^{q}\star_N \mathcal{X}^{1/2},
\end{eqnarray}
where $q \in \mathbb{R}$.

\begin{example}\label{exp:weighted power mean}
Given a probability vector $\bm{w}=[w_1,\cdots,w_k]$ with length $k$, $k$ random tensors $\mathcal{A}_i \in S_{\mbox{\tiny PD}}$ for $1 \leq i \leq k$ and $q \in [-1,1]\backslash \{0\}$, the weighted power mean for $0 < q  \leq 1$, represented by $\mathfrak{P}_{\bm{w},q}(\mathcal{A}_1,\cdots,\mathcal{A}_k)$, is expressed by
\begin{eqnarray}\label{eq1:exp:weighted power mean}
\mathcal{X}&=&\mathfrak{M}_{A,\bm{w}}(\mathcal{X}\#_q \mathcal{A}_1,\cdots,\mathcal{X}\#_q \mathcal{A}_k);
\end{eqnarray}
and,  for $-1 \leq q  < 0$, $\mathfrak{P}_{\bm{w},q}(\mathcal{A}_1,\cdots,\mathcal{A}_k)$ is expressed by
\begin{eqnarray}\label{eq2:exp:weighted power mean}
\mathcal{X}&=&\mathfrak{M}_{H,\bm{w}}(\mathcal{X}\#_{-q} \mathcal{A}_1,\cdots,\mathcal{X}\#_{-q} \mathcal{A}_k).
\end{eqnarray}
For $0 < q  \leq 1$, we have $\mathfrak{P}_{\bm{w},-q}(\mathcal{A}_1,\cdots,\mathcal{A}_k)=\mathfrak{P}^{*}_{\bm{w},q}(\mathcal{A}_1,\cdots,\mathcal{A}_k)$. The multivariate mean $\mathfrak{P}_{\bm{w},q}(\mathcal{A}_1,\cdots,\mathcal{A}_k)$ satisfies the multivariate Kubo-Ando axioms given at Section~\eqref{list:Kubo-Ando axioms mul Arg}. 
\end{example}

\begin{example}\label{exp:Karcher mean}
Given a probability vector $\bm{w}=[w_1,\cdots,w_k]$ with length $k$, $k$ random tensors $\mathcal{A}_i \in S_{\mbox{\tiny PD}}$ for $1 \leq i \leq k$, the Karcher equation is expressed by
\begin{eqnarray}\label{eq1:exp:Karcher mean}
\sum\limits_{i=1}^{k}w_i \log(\mathcal{X}^{-1/2}\star_N \mathcal{A}_i \star_N \mathcal{X}^{-1/2})
=\mathcal{O}.
\end{eqnarray}
The unique solution to the Karcher equation above is called the Karcher mean, represented by $\mathfrak{G}_{\bm{w}}(\mathcal{A}_1,\cdots,\mathcal{A}_k)$~\cite{lawson2014karcher}. In~\cite{lawson2014karcher}, they show that the multivariate mean $\mathfrak{G}_{\bm{w}}(\mathcal{A}_1,\cdots,\mathcal{A}_k)$ satisfies the multivariate Kubo-Ando axioms given at Section~\eqref{list:Kubo-Ando axioms mul Arg}. Besides, they also show the following inequalities:
\begin{eqnarray}
\mathfrak{P}_{\bm{w},-q}(\mathcal{A}_1,\cdots,\mathcal{A}_k) \preceq \mathfrak{G}_{\bm{w}}(\mathcal{A}_1,\cdots,\mathcal{A}_k) \preceq \mathfrak{P}_{\bm{w},q}(\mathcal{A}_1,\cdots,\mathcal{A}_k),
\end{eqnarray}
where $0 < q \leq 1$.
\end{example}

\section{Continuity and Differentiability of Karcher Means}\label{sec:Continuity and Differentiability of Karcher Means}

Theorem 6.3 in~\cite{lawson2014karcher} shows that the Karcher mean is continuous and differentiable in some neighborhood of the diagonal format. All tensors are with dimensions as $\mathbb{C}^{I_1 \times \dots \times I_N \times I_1 \times \dots \times I_N}$ in this section. Given a probability vector $\bm{w}=[w_1,\cdots,w_k]$ with length $k$ and $k$ identity tensors $\mathcal{I}$, the following map
\begin{eqnarray}\label{eq:diagonal map}
(\underbrace{\mathcal{I},\cdots,\mathcal{I}}_{k~\mbox{terms}})\rightarrow 
\mathfrak{G}_{\bm{w}}(\underbrace{\mathcal{I},\cdots,\mathcal{I}}_{k~\mbox{terms}}),
\end{eqnarray} 
is continuous and differentiable in some neighborhood of $(\underbrace{\mathcal{I},\cdots,\mathcal{I}}_{k~\mbox{terms}})$. In this section, we wish to extend this diagonal format to the following general format, i.e., the following map
\begin{eqnarray}\label{eq:general map}
(\mathcal{A}_1,\cdots,\mathcal{A}_{k})\rightarrow 
\mathfrak{G}_{\bm{w}}(\mathcal{A}_1,\cdots,\mathcal{A}_{k}),
\end{eqnarray} 
is continuous and differentiable in some neighborhood of any $(\mathcal{A}_1,\cdots,\mathcal{A}_{k}) \in S^k_{\mbox{\tiny PD}}$ in the sense of RTT. The main technique used in this section is based on \emph{prederiivative} of the function introduced in~\cite{pales1997inverse}. 

Let $\mathrm{X}, \mathrm{Y}$ be two Hilbert spaces, and $\mathfrak{A}$ be a collection of all random linear tensors (operators) from the space $\mathrm{X}$ to the space $\mathrm{Y}$. The  noncompactness measure of $\mathfrak{A}$, denoted as $\beta(\mathfrak{A})$, is defind by    
\begin{eqnarray}\label{eq:noncompact measure def}
\beta(\mathfrak{A}) \define \inf\{r| \exists k \in \mathbb{N}, \exists \mathcal{A}_1, \cdots, \mathcal{A}_k \in \mathfrak{A}~\mbox{such that}~\mathfrak{A} \subset \bigcup\limits_{i=1}^k B(A_i, r)\}.  
\end{eqnarray}
If $\mathfrak{A}$ is a compact space, $\beta(\mathfrak{A}) = 0$. The other measure is about the surjectivity measure of $\mathcal{A}$, denoted as $\gamma(\mathfrak{A})$, which is defined by
\begin{eqnarray}\label{eq:surjectivity measure def}
\gamma(\mathcal{A}) \define \sup\{r \in \mathbb{R}~\mbox{such that}~ r \mathbb{B}_{\mathrm{Y}} \subset \mathcal{A}(\mathbb{B}_{\mathrm{X}})\},
\end{eqnarray}
where $\mathcal{A} \in \mathfrak{A}$, and $\mathbb{B}_{\mathrm{X}}, \mathbb{B}_{\mathrm{Y}}$ are any unit norm balls inside spaces of $\mathrm{X}, \mathrm{Y}$, respectively. Then, $\gamma(\mathfrak{A})\define \inf\limits_{\mathcal{A} \in \mathfrak{A}}\gamma(\mathcal{A})$. 

Given two PD tensors $\mathcal{X}'$ and $\mathcal{X}''$, we say that these two tensors are \emph{square root representable} if we have 
\begin{eqnarray}\label{eq:square root representable}
\mathcal{X}'-\mathcal{X}'' &=& (\mathcal{X}'^{1/2}-\mathcal{X}''^{1/2}) (\mathcal{X}'^{1/2}+\mathcal{X}''^{1/2} + \mathcal{C}_{\mathcal{X}',\mathcal{X}''}),  
\end{eqnarray}
where $\mathcal{C}_{\mathcal{X}',\mathcal{X}''} \in \mathfrak{A}$. Note that the tensor $\mathcal{C}_{\mathcal{X}',\mathcal{X}''} \in \mathfrak{A} = \mathcal{O}$ if $\mathcal{X}' \star_N \mathcal{X}'' = \mathcal{X}'' \star_N \mathcal{X}'$, i.e., commutative. 

\begin{lemma}\label{lma:Thm 3, req iii}
Let $D \subset S_{\mbox{\tiny PD}}$ be an open set that includes $\mathfrak{G}_{\bm{w}}(\hat{\mathcal{A}}_1,\cdots,\hat{\mathcal{A}}_k) \in D$ and all pair of distinct tensors within $D$ are square root representable defined by Eq.~\eqref{eq:square root representable}, where $\mathfrak{G}_{\bm{w}}(\hat{\mathcal{A}}_1,\cdots,\hat{\mathcal{A}}_k)$ is the Karcher mean with respect to the weight $\bm{w}$ and $(\hat{\mathcal{A}}_1,\cdots,\hat{\mathcal{A}}_k) \in S^k_{\mbox{\tiny PD}}$.  We define the following tensor-valued function:
\begin{eqnarray}\label{eq1:lma:Thm 3, req iii}
F_{\bm{w}}(\mathcal{A}_1,\cdots,\mathcal{A}_{k},\mathcal{X})\define\sum\limits_{i=1}^k w_i \log(\mathcal{X}^{1/2}\star_N \mathcal{A}^{-1}_i \star_N \mathcal{X}^{1/2}).
\end{eqnarray}
Then, we have the following property with respect to the function $F_{\bm{w}}(\mathcal{A}_1,\cdots,\mathcal{A}_{k},\mathcal{X})$. For all $\epsilon > 0$, there exists a neighborhood $W_{\epsilon} \subset  S^k_{\mbox{\tiny PD}} \times D$ of $(\hat{\mathcal{A}}_1,\cdots,\hat{\mathcal{A}}_k,\mathfrak{G}_{\bm{w}}(\hat{\mathcal{A}}_1,\cdots,\hat{\mathcal{A}}_k))$ such that, for any $(\mathcal{A}_1,\cdots,\mathcal{A}_{k},\mathcal{X}'), (\mathcal{A}_1,\cdots,\mathcal{A}_{k},\mathcal{X}'') \in W_{\epsilon}$, we have:  
\begin{eqnarray}\label{eq2:lma:Thm 3, req iii}
F_{\bm{w}}(\mathcal{A}_1,\cdots,\mathcal{A}_{k},\mathcal{X}') -F_{\bm{w}}(\mathcal{A}_1,\cdots,\mathcal{A}_{k},\mathcal{X}'') &=& \acute{\mathcal{A}} \star_N (\mathcal{X}' - \mathcal{X}'') + \epsilon \left\Vert  \mathcal{X}' - \mathcal{X}'' \right\Vert_{\rho}\mathcal{B}
\end{eqnarray}
where $\acute{\mathcal{A}}, \mathcal{B} \in \mathfrak{A}$, and $\mathcal{B}$ is a tensor with norm $\left\Vert \mathcal{B} \right\Vert_{\rho} \leq 1$.
\end{lemma}
\textbf{Proof:}
In the region $D$ around $\mathfrak{G}_{\bm{w}}(\hat{\mathcal{A}}_1,\cdots,\hat{\mathcal{A}}_k)$, we select a neighborhood $D_\delta$ of $\mathcal{X}', \mathcal{X}''$ such that 
\begin{eqnarray}\label{eq3:lma:Thm 3, req iii}
\mathcal{X}'^{1/2} = \mathcal{X}''^{1/2} + \Delta,
\end{eqnarray}
where the tensor $\Delta$ satisfies with $\left\Vert \Delta \right\Vert_{\rho} < \delta$ with $\delta$ having 
\begin{eqnarray}\label{eq3.1:lma:Thm 3, req iii}
\delta^2 \leq C\epsilon \left\Vert \mathcal{X}' - \mathcal{X}'' \right\Vert_{\rho}, 
\end{eqnarray}
where $C>0$ is a constant indepdent of selection $\mathcal{X}'$ and $\mathcal{X}''$ in $D$.
 
 and $\delta$. From the definition provided by Eq.~\eqref{eq1:lma:Thm 3, req iii}, we have
\begin{eqnarray}\label{eq4:lma:Thm 3, req iii}
\lefteqn{F_{\bm{w}}(\mathcal{A}_1,\cdots,\mathcal{A}_{k},\mathcal{X}') -F_{\bm{w}}(\mathcal{A}_1,\cdots,\mathcal{A}_{k},\mathcal{X}'')}\nonumber \\
&=&\sum\limits_{i=1}^k w_i \log\frac{(\mathcal{X}''^{1/2} + \Delta)\mathcal{A}^{-1}_i (\mathcal{X}''^{1/2} + \Delta)}{\mathcal{X}''^{-1/2}\mathcal{A}^{-1}_i \mathcal{X}''^{-1/2}}\nonumber \\
&=& \sum\limits_{i=1}^k w_i \log\left(\mathcal{I}+\frac{\Delta \mathcal{A}^{-1}_i\mathcal{X}''^{1/2}}{\mathcal{X}''^{1/2}\mathcal{A}_i^{-1} \mathcal{X}''^{1/2}}+ \frac{\mathcal{X}''^{1/2}\mathcal{A}^{-1}_i \Delta}{\mathcal{X}''^{1/2}\mathcal{A}_i^{-1} \mathcal{X}''^{1/2}}+\frac{\Delta \mathcal{A}^{-1}_i \Delta}{\mathcal{X}''^{1/2}\mathcal{A}_i^{-1} \mathcal{X}''^{1/2}}\right)\nonumber \\
&=_1&\sum\limits_{i=1}^k w_i  \left(\widetilde{\frac{ \mathcal{A}^{-1}_i\mathcal{X}''^{1/2}}{\mathcal{X}''^{1/2}\mathcal{A}_i^{-1} \mathcal{X}''^{1/2}}}+ \frac{\mathcal{X}''^{1/2}\mathcal{A}^{-1}_i }{\mathcal{X}''^{1/2}\mathcal{A}_i^{-1} \mathcal{X}''^{1/2}}\right)\Delta + \mathrm{O}(\Delta^2)\nonumber \\
&=_2&\sum\limits_{i=1}^k w_i \left(\widetilde{\frac{ \mathcal{A}^{-1}_i\mathcal{X}''^{1/2}}{\mathcal{X}''^{1/2}\mathcal{A}_i^{-1} \mathcal{X}''^{1/2}}}+ \frac{\mathcal{X}''^{1/2}\mathcal{A}^{-1}_i }{\mathcal{X}''^{1/2}\mathcal{A}_i^{-1} \mathcal{X}''^{1/2}}\right)\frac{(\mathcal{X}' - \mathcal{X}'')}{(\mathcal{X}'^{1/2} + \mathcal{X}''^{1/2}+ \mathcal{C}_{\mathcal{X}',\mathcal{X}''})}\nonumber \\
&& + \mathrm{O}(\Delta^2)\nonumber \\
&=_3&\acute{\mathcal{A}}(\mathcal{X}' - \mathcal{X}'')+\epsilon \left\Vert  \mathcal{X}' - \mathcal{X}'' \right\Vert_{\rho}\mathcal{B}
\end{eqnarray}
where $=_1$ comes from the Taylor expansion for $\log(1+x)=\sum\limits_{i=1}^{\infty}(-1)^{i-1}\frac{x^i}{i}$ and $\frac{\Delta \mathcal{A}^{-1}_i\mathcal{X}''^{1/2}}{\mathcal{X}''^{1/2}\mathcal{A}_i^{-1} \mathcal{X}''^{1/2}} = \widetilde{\frac{ \mathcal{A}^{-1}_i\mathcal{X}''^{1/2}}{\mathcal{X}''^{1/2}\mathcal{A}_i^{-1} \mathcal{X}''^{1/2}}}\Delta$, $=_2$ uses the fact that tensors $\mathcal{X}'$ and $\mathcal{X}''$ are \emph{square root representable}, and $=_3$ comes from Eq.~\eqref{eq3.1:lma:Thm 3, req iii} and summands rearrangments by setting 
\begin{eqnarray}
\acute{\mathcal{A}} \define \left[\sum\limits_{i=1}^k w_i \left(\widetilde{\frac{ \mathcal{A}^{-1}_i\mathcal{X}''^{1/2}}{\mathcal{X}''^{1/2}\mathcal{A}_i^{-1} \mathcal{X}''^{1/2}}}+ \frac{\mathcal{X}''^{1/2}\mathcal{A}^{-1}_i }{\mathcal{X}''^{1/2}\mathcal{A}_i^{-1} \mathcal{X}''^{1/2}}\right)\right]\left(\mathcal{X}'^{1/2} + \mathcal{X}''^{1/2}+ \mathcal{C}_{\mathcal{X}',\mathcal{X}''}\right)^{-1}.
\end{eqnarray}
Then, this Lemma is proved by selecting  $W_{\epsilon} =  \tilde{S}_{\mbox{\tiny PD,1}} \times \cdots \times \tilde{S}_{\mbox{\tiny PD,k}} \times D_\delta$, where $\tilde{S}_{\mbox{\tiny PD,i}}$ is the open set in $S_{\mbox{\tiny PD}}$ that contains $\hat{\mathcal{A}}_i$ for $1\leq i \leq k$.  
$\hfill \Box$

Following Lemma is about the implicit function of $F_{\bm{w}}(\mathcal{A}_1,\cdots,\mathcal{A}_{k},\mathcal{X})$ if Lemma~\ref{lma:Thm 3, req iii} is satisfied. We define $\mathfrak{C}(\mathrm{X},\mathrm{Y})$ as a collection of all continuous functions mapping from the space $\mathrm{X}$ to the space $\mathrm{Y}$.

\begin{lemma}\label{lma:Thm 3}
Let $D \subset S_{\mbox{\tiny PD}}$ be an open set that includes $\mathfrak{G}_{\bm{w}}(\hat{\mathcal{A}}_1,\cdots,\hat{\mathcal{A}}_k) \in D$ and all pair of distinct tensors within $D$ are square root representable defined by Eq.~\eqref{eq:square root representable}, where $\mathfrak{G}_{\bm{w}}(\hat{\mathcal{A}}_1,\cdots,\hat{\mathcal{A}}_k)$ is the Karcher mean with respect to the weight $\bm{w}$ and $(\hat{\mathcal{A}}_1,\cdots,\hat{\mathcal{A}}_k) \in S^k_{\mbox{\tiny PD}}$.  We define the following tensor-valued function:
\begin{eqnarray}\label{eq1:lma:Thm 3}
F_{\bm{w}}(\mathcal{A}_1,\cdots,\mathcal{A}_{k},\mathcal{X})\define\sum\limits_{i=1}^k w_i \log(\mathcal{X}^{1/2}\star_N \mathcal{A}^{-1}_i \star_N \mathcal{X}^{1/2}).
\end{eqnarray}
For notation simplicity, we set $\underline{\hat{\mathcal{A}}} \define (\hat{\mathcal{A}}_1,\cdots,\hat{\mathcal{A}}_k)$. Then, for all $\beta(\mathfrak{A}) < \tau_1 < \tau_2 < \gamma(\mathfrak{A})$, there exists a neighborhood $U$ of $(\underline{\hat{\mathcal{A}}},\mathfrak{G}_{\bm{w}}(\underline{\hat{\mathcal{A}}}),F_{\bm{w}}(\underline{\hat{\mathcal{A}}},\mathfrak{G}_{\bm{w}}(\underline{\hat{\mathcal{A}}})))$ and a function $\phi:U \rightarrow D$ such that 
\begin{eqnarray}\label{eq:10-1}
F_{\bm{w}}(\underline{\mathcal{A}},\phi(\underline{\mathcal{A}},\mathcal{X},\mathcal{Y}))&=&\mathcal{Y},
\end{eqnarray}
and
\begin{eqnarray}\label{eq:10-2}
\left\Vert \mathcal{X} -  \phi(\underline{\mathcal{A}},\mathcal{X},\mathcal{Y}) \right\Vert_{\rho}
&\leq& \frac{\left\Vert F_{\bm{w}}(\underline{\mathcal{A}},\mathcal{X}) - \mathcal{Y}\right\Vert_{\rho}}{\tau_2 - \tau_1},
\end{eqnarray}
are satisfied by all $(\underline{\mathcal{A}},\mathcal{X},\mathcal{Y}) \in U$.
\end{lemma}
\textbf{Proof:}
Pick a positive number $r$ betwee $\beta(\mathfrak{A})$ and $\tau_1$ and the definition of $\beta(\mathfrak{A})$ provied by Eq.~\eqref{eq:noncompact measure def}, we have
\begin{eqnarray}\label{eq4:lma:Thm 3}
\mathfrak{A} \subset \bigcup\limits_{i=1}^k B(\mathcal{A}_i,r).
\end{eqnarray}
Let $\Upsilon_{\underline{\mathcal{A}}}$ be the convex hull of $\{\mathcal{A}_1, \cdots, \mathcal{A}_k\}$. Then $\Upsilon_{\underline{\mathcal{A}}} \in \mathfrak{A}$ is compact and convex. By setting $\epsilon_1 = \tau_1 - r$, from
~Lemma~\ref{lma:Thm 3, req iii}, we have a neighborhood $W_{\epsilon_1} \subset S^k_{\mbox{\tiny PD}} \times D$ such that 
\begin{eqnarray}\label{eq5:lma:Thm 3}
F_{\bm{w}}(\underline{\mathcal{A}},\mathcal{X}')-F_{\bm{w}}(\underline{\mathcal{A}},\mathcal{X}'')
&=&  \acute{\mathcal{A}} (\mathcal{X}' - \mathcal{X}'')+ \epsilon_1 \left\Vert  \mathcal{X}' - \mathcal{X}'' \right\Vert_{\rho}\mathcal{B},
\end{eqnarray}
where $(\underline{\mathcal{A}},\mathcal{X}'), (\underline{\mathcal{A}},\mathcal{X}'') \in W_{\epsilon_1}$. Note that $\mathfrak{A} \subset \Upsilon_{\underline{\mathcal{A}}} + r \mathbb{B}_{L(S_{\mbox{\tiny PD}},S_{\mbox{\tiny PD}})}$, where $ \mathbb{B}_{L(S_{\mbox{\tiny PD}},S_{\mbox{\tiny PD}})}$  is the set of all tensors (opertors) from $S_{\mbox{\tiny PD}}$ to $S_{\mbox{\tiny PD}}$ with the norm measured by $\left\Vert \cdot \right\Vert_{\rho}$ less or equal then one. 

Pick any $\upsilon \in \Upsilon_{\underline{\mathcal{A}}}$, since the mapping $\upsilon \rightarrow \gamma(\upsilon)$ is continuous on $\Upsilon_{\underline{\mathcal{A}}}$ and the set $\Upsilon_{\underline{\mathcal{A}}}$ is compact, we can find a constant $c$ such that $\tau_2 < c < \gamma(\upsilon)$ for all $\upsilon \in \Upsilon_{\underline{\mathcal{A}}}$. We can have a surjective map $g: \mathbb{C}(\Upsilon_{\underline{\mathcal{A}}},S_{\mbox{\tiny PD}}) \rightarrow \mathbb{C}(\Upsilon_{\underline{\mathcal{A}}},S_{\mbox{\tiny PD}})$ such that 
\begin{eqnarray}\label{eq:12}
c \mathbb{B}_{\mathbb{C}(\Upsilon_{\underline{\mathcal{A}}},S_{\mbox{\tiny PD}})} \subset g(\mathbb{B}_{\mathbb{C}(\Upsilon_{\underline{\mathcal{A}}},S_{\mbox{\tiny PD}})}).
\end{eqnarray}
where $0 < c < \gamma(g)$. 

Set $\epsilon_2 = \tau_2 - \tau_1$ and select $\delta > 0$ with a neighborhood $U$ of $(\underline{\hat{\mathcal{A}}},\mathfrak{G}_{\bm{w}}(\underline{\hat{\mathcal{A}}}),F_{\bm{w}}(\underline{\hat{\mathcal{A}}},\mathfrak{G}_{\bm{w}}(\underline{\hat{\mathcal{A}}})))$ such that 
\begin{eqnarray}\label{eq:12}
(\underline{\mathcal{A}},\mathcal{T}) \in W_{\epsilon_1},~\left\Vert \mathcal{X} - \mathfrak{G}_{\bm{w}}(\underline{\hat{\mathcal{A}}})\right\Vert_{\rho} \leq \frac{\delta}{3},~\left\Vert F_{\bm{w}}(\underline{\mathcal{A}},\mathcal{X}) - \mathcal{Y}\right\Vert_{\rho}\leq \frac{\epsilon_2 \delta}{3},
\end{eqnarray}
where $(\underline{\mathcal{A}}, \mathcal{X}, \mathcal{Y}) \in U$ and $\mathcal{T} \in \mathfrak{G}_{\bm{w}}(\underline{\hat{\mathcal{A}}}) + \delta \mathbb{B}_{S_{\mbox{\tiny PD}}}$.

We will apply Ekeland's variational principle to obtain the function $\phi$ that satisfies Eqs.~\eqref{eq:10-1} and~\eqref{eq:10-2}~\cite{ekeland1974variational}. By applying Ekeland's variational principle to the complete metric space $\mathfrak{G}_{\bm{w}}(\underline{\hat{\mathcal{A}}}) + \delta \mathbb{B}_{S_{\mbox{\tiny PD}}}$ and the mapping $\mathcal{T} \rightarrow \left\Vert F_{\bm{w}}(\underline{\mathcal{A}},\mathcal{T}) - \mathcal{Y}\right\Vert_{\rho}$, there exists 
the function $\phi \in \mathfrak{G}_{\bm{w}}(\underline{\hat{\mathcal{A}}}) + \delta \mathbb{B}_{S_{\mbox{\tiny PD}}}$ satisfying the following relations:
\begin{eqnarray}\label{eq:13}
\left\Vert F_{\bm{w}}(\underline{\mathcal{A}},\phi) - \mathcal{Y} \right\Vert_{\rho}
+\epsilon_2 \left\Vert \mathcal{X} - \phi \right\Vert_{\rho} &\leq& \left\Vert F_{\bm{w}}(\underline{\mathcal{A}},\mathcal{X}) - \mathcal{Y}\right\Vert_{\rho},
\end{eqnarray}
and
\begin{eqnarray}\label{eq:14}
\left\Vert F_{\bm{w}}(\underline{\mathcal{A}},\phi) - \mathcal{Y} \right\Vert_{\rho}&\leq&
\left\Vert F_{\bm{w}}(\underline{\mathcal{A}},\mathcal{T}) - \mathcal{Y} \right\Vert_{\rho} + \epsilon_2 \left\Vert \mathcal{T} - \phi\right\Vert_{\rho},
\end{eqnarray}
where $\mathcal{T} \in \mathfrak{G}_{\bm{w}}(\underline{\hat{\mathcal{A}}}) + \delta \mathbb{B}_{S_{\mbox{\tiny PD}}}$. Then, the inequality provided by Eq.~\eqref{eq:13} implies Eq.~\eqref{eq:10-2}. Now, we want to prove  Eq.~\eqref{eq:10-1}. 

From Eq.~\eqref{eq:13}, we have
\begin{eqnarray}
\left\Vert \mathcal{X} - \phi \right\Vert_{\rho} \leq \frac{\left\Vert F_{\bm{w}}(\underline{\mathcal{A}},\mathcal{X}) - \mathcal{Y}\right\Vert_{\rho}}{\epsilon_2} \leq \frac{\delta}{3},
\end{eqnarray}
combined with the condition $\left\Vert \mathcal{X} - \mathfrak{G}_{\bm{w}}(\underline{\hat{\mathcal{A}}})\right\Vert_{\rho} \leq \frac{\delta}{3}$, we obtain 
\begin{eqnarray}\label{eq:14-1}
\left\Vert  \phi -\mathfrak{G}_{\bm{w}}(\underline{\hat{\mathcal{A}}}) \right\Vert_{\rho} \leq \frac{2\delta}{3}.
\end{eqnarray}

From Eq.~\eqref{eq:12}, there exists an element $\mathcal{X}_{\upsilon}$ for all $\epsilon \in \Upsilon_{\underline{\mathcal{A}}}$ such that the map $\upsilon \rightarrow \mathcal{X}_{\upsilon}$ is continuous and the following relations: 
\begin{eqnarray}\label{eq:15-1}
\mathcal{Y} - F_{\bm{w}}(\underline{\mathcal{A}},\phi) &=& \upsilon \mathcal{X}_{\upsilon}, 
\end{eqnarray}
and
\begin{eqnarray}\label{eq:15-2}
\left\Vert \mathcal{X}_{\upsilon}\right\Vert_{\rho} &<& \frac{\delta}{3}, 
\end{eqnarray}
where $\upsilon \in \Upsilon_{\underline{\mathcal{A}}}$. If we set $\mathcal{T} = \phi + \mathcal{X}_{\upsilon}$, from Eq.~\eqref{eq:14-1} and Eq.~\eqref{eq:15-2}, we have $\mathcal{T} \in \mathfrak{G}_{\bm{w}}(\underline{\hat{\mathcal{A}}}) + \delta \mathbb{B}_{S_{\mbox{\tiny PD}}}$, then, from Eq.~\eqref{eq:14}, we have 
\begin{eqnarray}\label{eq:16}
\left\Vert F_{\bm{w}}(\underline{\mathcal{A}},\phi) - \mathcal{Y}\right\Vert_{\rho}
&\leq& \left\Vert F_{\bm{w}}(\underline{\mathcal{A}},\phi + \mathcal{X}_{\upsilon}) -  F_{\bm{w}}(\underline{\mathcal{A}},\phi) - \upsilon \mathcal{X}_{\upsilon}\right\Vert_{\rho} + \epsilon_2 \left\Vert \mathcal{X}_{\upsilon} \right\Vert_{\rho},
\end{eqnarray}
where $\upsilon \in \Upsilon_{\underline{\mathcal{A}}}$.

We define the set-valued map $\Phi:\Upsilon_{\underline{\mathcal{A}}} \rightarrow 2^{\Upsilon_{\underline{\mathcal{A}}}}$ by 
\begin{eqnarray}
\Phi(\upsilon) \define \{\varrho \in \Upsilon_{\underline{\mathcal{A}}}| \left\Vert F_{\bm{w}}(\underline{\mathcal{A}},\phi + \mathcal{X}_{\upsilon}) -  F_{\bm{w}}(\underline{\mathcal{A}},\phi) - \varrho \mathcal{X}_{\upsilon}\right\Vert_{\rho}  \leq \tau_1 \left\Vert\mathcal{X}_{\upsilon}\right\Vert_{\rho}\}.
\end{eqnarray}
From Lemma~\ref{lma:Thm 3, req iii}, we know that the set $\Phi(\upsilon)$ is nonempty for all $\upsilon \in \Upsilon_{\underline{\mathcal{A}}}$. Besides, the set $\Phi(\upsilon)$ is a convex and compact set. Before applying Kakutani’s fixed point theorem to the set $\Phi(\upsilon)$, we have to show that the set $\Phi(\upsilon)$ is semicontinuous. Let $\upsilon_0 \in \Upsilon_{\underline{\mathcal{A}}}$ be an arbitrary element, $\upsilon_n, \varrho_n \in \Phi(\upsilon_n)$ are two sequences with $\upsilon_0 = \lim\limits_{n \rightarrow \infty} \upsilon_n$ and $\varrho_0 = \lim\limits_{n \rightarrow \infty} \varrho_n$. Since we have
\begin{eqnarray}
\left\Vert F_{\bm{w}}(\underline{\mathcal{A}},\phi + \mathcal{X}_{\upsilon_n}) - F_{\bm{w}}(\underline{\mathcal{A}},\phi) - \varrho_n \mathcal{X}_{\upsilon_n}\right\Vert_{\rho} \leq \tau_1 \left\Vert \mathcal{X}_{\upsilon_n}\right\Vert_{\rho},
\end{eqnarray}
by taking $n \rightarrow \infty$ and using the continuity of $\upsilon \rightarrow \mathcal{X}_{\upsilon}$, we have $\varrho_0 \in \Phi(\upsilon_0)$. This shows that the set $\Phi(\upsilon)$ is semicontinuous. 

From Kakutani’s fixed point theorem, there exists a fixed element $\tilde{\upsilon} \in \Phi(\upsilon)$, i.e., 
\begin{eqnarray}
\left\Vert F_{\bm{w}}(\underline{\mathcal{A}},\phi + \mathcal{X}_{\tilde{\upsilon}}) - F_{\bm{w}}(\underline{\mathcal{A}},\phi) - \tilde{\upsilon}\mathcal{X}_{\tilde{\upsilon}}\right\Vert_{\rho} \leq \tau_1 \left\Vert \mathcal{X}_{\tilde{\upsilon}}\right\Vert_{\rho}.
\end{eqnarray}
By setting $\upsilon = \tilde{\upsilon}$ into Eq.~\eqref{eq:16} and applying the inequality provided by Eq.~\eqref{eq:15-2}, we get
\begin{eqnarray}\label{eq-last:lma:Thm 3}
\left\Vert F_{\bm{w}}(\underline{\mathcal{A}},\phi) - \mathcal{Y}\right\Vert_{\rho}
\leq (\tau_1+\epsilon_2)\left\Vert\mathcal{X}_{\tilde{\upsilon}}\right\Vert_{\rho} = \tau_2 \left\Vert\mathcal{X}_{\tilde{\upsilon}}\right\Vert_{\rho} 
\leq \frac{\tau_2 \left\Vert F_{\bm{w}}(\underline{\mathcal{A}},\phi) - \mathcal{Y}\right\Vert_{\rho}}{c}.
\end{eqnarray}
Because $\tau_2 < c$, the inequality given by Eq.~\eqref{eq-last:lma:Thm 3} yields $\left\Vert F_{\bm{w}}(\underline{\mathcal{A}},\phi) - \mathcal{Y}\right\Vert_{\rho}=0$. This shows the existence of $\phi$ that satifies Eq.~\eqref{eq:10-1}. 
$\hfill \Box$

Before presenting our main theorem in this section about the continuity and the differentiability of Karcher Means, we need to present one more lemma about the surjective measure $\gamma$.
\begin{lemma}\label{lma:4}
Let $\mathrm{X}, \mathrm{Y}$ be two Hilbert spaces of tensors and $\mathcal{A}: \mathrm{X}\rightarrow\mathrm{Y}$ be an injective bouded tensor, we have
\begin{eqnarray}\label{eq1:lma:4}
\left\Vert \mathcal{A}\mathcal{X} \right\Vert_{\rho} \geq \gamma(\mathcal{A})\left\Vert \mathcal{X} \right\Vert_{\rho},
\end{eqnarray}
where $\forall \mathcal{X} \in \mathrm{X}$.
\end{lemma}
\textbf{Proof:}
If $\gamma(\mathcal{A})=0$, this makes Eq.~\eqref{eq1:lma:4} valid. We will prove this Lemma by contradiction. Suppose Eq.~\eqref{eq1:lma:4} is not valid for $\gamma(\mathcal{A}) > 0$, then there exists a positive number $c <  \gamma(\mathcal{A})$ and an $\mathcal{X} \in \mathrm{X}$ such that 
\begin{eqnarray}\label{eq2:lma:4}
\left\Vert \mathcal{A}\mathcal{X} \right\Vert_{\rho} < c\left\Vert \mathcal{X} \right\Vert_{\rho}.
\end{eqnarray}
Then, we have
\begin{eqnarray}\label{eq3:lma:4}
\frac{c}{\left\Vert \mathcal{A} \right\Vert_{\rho}} \mathcal{A}\mathcal{X} \in c \mathbb{B}_{\mathrm{Y}}
\subset \mathcal{A}\mathbb{B}_{\mathrm{X}}, 
\end{eqnarray}
which is equivalent to have $\tilde{\mathcal{X}}$ such that 
\begin{eqnarray}\label{eq4:lma:4}
\mathcal{A}\left(\frac{c}{\left\Vert \mathcal{A}\mathcal{X}\right\Vert_{\rho}} \mathcal{X}\right)=\mathcal{A}(\tilde{\mathcal{X}}).
\end{eqnarray}
From the injectivity of $\mathcal{A}$, we have 
\begin{eqnarray}\label{eq5:lma:4}
\frac{c}{\left\Vert \mathcal{A}\mathcal{X}\right\Vert_{\rho}} \mathcal{X}&=&\tilde{\mathcal{X}}.
\end{eqnarray}
By taking the $\left\Vert \cdot \right\Vert_{\rho}$ operation at both sides of Eq.~\eqref{eq5:lma:4}, we have
\begin{eqnarray}\label{eq5:lma:4}
1 < \frac{c\left\Vert \mathcal{X}\right\Vert_{\rho}}{\left\Vert \mathcal{A}\mathcal{X}\right\Vert_{\rho}}&=&\left\Vert\tilde{\mathcal{X}}\right\Vert_{\rho} \leq 1.
\end{eqnarray}
However, this contradics to our assumption. 
$\hfill \Box$

We are ready to present our main theorem about the continuity and the differentiability of Karcher Means in this section. For notational simplicity, we set $\underline{\hat{\mathcal{A}}} \define (\hat{\mathcal{A}}_1,\cdots,\hat{\mathcal{A}}_k)$. 
\begin{theorem}\label{thm4:Continuity and Differentiability of Karcher Means}
Let $D \subset S_{\mbox{\tiny PD}}$ be an open set that includes $\mathfrak{G}_{\bm{w}}(\underline{\hat{\mathcal{A}}}) \in D$ and all pair of distinct tensors within $D$ are square root representable defined by Eq.~\eqref{eq:square root representable}, where $\mathfrak{G}_{\bm{w}}(\underline{\hat{\mathcal{A}}})$ is the Karcher mean with respect to the weight $\bm{w}$ and $\underline{\hat{\mathcal{A}}} \in S^k_{\mbox{\tiny PD}}$.  We define the following tensor-valued function:
\begin{eqnarray}\label{eq1:Thm4}
F_{\bm{w}}(\underline{\mathcal{A}},\mathcal{X})\define\sum\limits_{i=1}^k w_i \log(\mathcal{X}^{1/2}\star_N \mathcal{A}^{-1}_i \star_N \mathcal{X}^{1/2}).
\end{eqnarray}
We assume that all tensors from $\mathfrak{A}$, i.e., a collection of all random linear tensors from the space $S_{\mbox{\tiny PD}}$ to the space $S_{\mbox{\tiny PD}}$, are injective. Then, there exist a neighborhood $V$ of $(\underline{\hat{\mathcal{A}}},F_{\bm{w}}(\underline{\hat{\mathcal{A}}},\mathfrak{G}_{\bm{w}}(\underline{\hat{\mathcal{A}}})))$ and a unique determined function $\psi: V \rightarrow D$ such that 
\begin{eqnarray}\label{eq2:Thm4 19}
F_{\bm{w}}(\underline{\mathcal{A}},\psi(\underline{\mathcal{A}},\mathcal{Y}))
= \mathcal{Y},
\end{eqnarray}
where $(\underline{\mathcal{A}},\mathcal{Y}) \in V$. The mapping below 
\begin{eqnarray}\label{eq3:Thm4 20}
\underline{\mathcal{A}} &\rightarrow&\psi(\underline{\mathcal{A}},F_{\bm{w}}(\underline{\hat{\mathcal{A}}},\mathfrak{G}_{\bm{w}}(\underline{\hat{\mathcal{A}}}))),
\end{eqnarray}
is continuous at $\underline{\mathcal{A}}=\underline{\hat{\mathcal{A}}}$. Finally, the mapping $\psi(\underline{\mathcal{A}},\mathcal{Y})$ is prederivative, i.e., for all $\epsilon > 0$, there exists a neighborhood $V_{\epsilon} \subset V$ such that 
\begin{eqnarray}\label{eq:Thm4 21}
\psi(\underline{\mathcal{A}},\mathcal{Y}') - \psi(\underline{\mathcal{A}},\mathcal{Y}'')
= \acute{\mathcal{A}}^{-1}\star_N (\mathcal{Y}' - \mathcal{Y}'')+ \epsilon \left\Vert \mathcal{Y}' - \mathcal{Y}'' \right\Vert_{\rho}\mathcal{B},
\end{eqnarray}
where $\acute{\mathcal{A}} \in \mathfrak{A}$ and $(\underline{\mathcal{A}},\mathcal{Y}'),(\underline{\mathcal{A}},\mathcal{Y}'') \in V_{\epsilon}$.
\end{theorem}
\textbf{Proof:}
Let $\beta(\mathfrak{A}) < r < \tau_1 < \tau_2 < \gamma(\mathfrak{A})$, $\epsilon_1=\tau_1 - r < \gamma(\mathfrak{A})$, and $W=W_{\epsilon_1}$ be the corresponding nieghborhood in Lemma~\ref{lma:Thm 3}. Our first goal is to show that if $F_{\bm{w}}(\underline{\mathcal{A}},\mathcal{X}') =F_{\bm{w}}(\underline{\mathcal{A}},\mathcal{X}'')$, this implies $\mathcal{X}' = \mathcal{X}''$. If $F_{\bm{w}}(\underline{\mathcal{A}},\mathcal{X}') =F_{\bm{w}}(\underline{\mathcal{A}},\mathcal{X}'')$, then, from Lemma~\ref{lma:Thm 3, req iii}, there exists $\mathcal{A} \in \mathfrak{A}$ and $\mathcal{Y}$ with $\left\Vert \mathcal{Y} \right\Vert_{\rho} \leq 1$ such that 
\begin{eqnarray}\label{eq5:Thm4}
\mathcal{A}(\mathcal{X}' - \mathcal{X}'') + \epsilon_1 \left\Vert \mathcal{X}' - \mathcal{X}'' \right\Vert_{\rho}\mathcal{Y}=\mathcal{O}.
\end{eqnarray}
From Lemma~\ref{lma:4}, we have 
\begin{eqnarray}
\epsilon_1 \left\Vert\mathcal{X}' - \mathcal{X}''\right\Vert_{\rho}&\geq&
\epsilon_1 \left\Vert \mathcal{Y} \right\Vert_{\rho}\left\Vert\mathcal{X}' - \mathcal{X}''\right\Vert_{\rho} \nonumber \\
&=& \left\Vert\mathcal{A}(\mathcal{X}' - \mathcal{X}'')\right\Vert_{\rho} \nonumber \\
&\geq&\gamma(\mathcal{A})\left\Vert \mathcal{X}' - \mathcal{X}'' \right\Vert_{\rho}.
\end{eqnarray}
If $\mathcal{X}' \neq \mathcal{X}''$,  we will obtain a contradiction as $\epsilon_1 \geq \gamma(\mathcal{A})$. Therefore,  if $F_{\bm{w}}(\underline{\mathcal{A}},\mathcal{X}') =F_{\bm{w}}(\underline{\mathcal{A}},\mathcal{X}'')$, we have $\mathcal{X}' = \mathcal{X}''$.

From Lemma~\ref{lma:Thm 3}, there exists a neighborhood $U$ and $\phi$ such that $(\underline{\mathcal{A}},\phi(\underline{\mathcal{A}},\mathcal{X},\mathcal{Y})) \in W_{\epsilon_1}$ given $(\underline{\mathcal{A}},\mathcal{X},\mathcal{Y}) \in U$. Define $V$ as the projection of $U$ on the space $S^k_{\mbox{\tiny PD}} \times S_{\mbox{\tiny PD}}$ and define $\psi: V \rightarrow D$ by  
\begin{eqnarray}
\psi(\underline{\mathcal{A}},\mathcal{Y}) \define \phi(\underline{\mathcal{A}},\mathcal{X},\mathcal{Y}),
\end{eqnarray}
where $(\underline{\mathcal{A}},\mathcal{X},\mathcal{Y}) \in U$. If $(\underline{\mathcal{A}},\mathcal{X}',\mathcal{Y})$ and $(\underline{\mathcal{A}},\mathcal{X}'',\mathcal{Y})$ are in $U$, then $(\underline{\mathcal{A}},\phi(\underline{\mathcal{A}},\mathcal{X}',\mathcal{Y}))$ and $(\underline{\mathcal{A}},\phi(\underline{\mathcal{A}},\mathcal{X}'',\mathcal{Y}))$ are in $W_{\epsilon_1}$. From Eq.~\eqref{eq:10-1}, we have
\begin{eqnarray}
F_{\bm{w}}(\underline{\mathcal{A}},\phi(\underline{\mathcal{A}},\mathcal{X}',\mathcal{Y}))=
\mathcal{Y} = F_{\bm{w}}(\underline{\mathcal{A}},\phi(\underline{\mathcal{A}},\mathcal{X}'',\mathcal{Y})).
\end{eqnarray}
We have Eq.~\eqref{eq2:Thm4 19} from $\phi(\underline{\mathcal{A}},\mathcal{X}',\mathcal{Y}) = \phi(\underline{\mathcal{A}},\mathcal{X}'',\mathcal{Y})$.

By setting $\mathcal{X} \define \mathfrak{G}_{\bm{w}}(\underline{\hat{\mathcal{A}}})$ and $\mathcal{Y} \define F_{\bm{w}}(\underline{\hat{\mathcal{A}}},\mathfrak{G}_{\bm{w}}(\underline{\hat{\mathcal{A}}}))$ into the inequality provided by Eq.~\eqref{eq:10-2}, we have
\begin{eqnarray}\label{eq22neg1:Thm4}
\left\Vert \mathfrak{G}_{\bm{w}}(\underline{\hat{\mathcal{A}}}) - \psi(\underline{\mathcal{A}}, 
F_{\bm{w}}(\underline{\hat{\mathcal{A}}},\mathfrak{G}_{\bm{w}}(\underline{\hat{\mathcal{A}}})))\right\Vert_{\rho} \leq \frac{\left\Vert F_{\bm{w}}(\underline{\mathcal{A}},\mathfrak{G}_{\bm{w}}(\underline{\hat{\mathcal{A}}}))-F_{\bm{w}}(\underline{\hat{\mathcal{A}}},\mathfrak{G}_{\bm{w}}(\underline{\hat{\mathcal{A}}}))\right\Vert_{\rho}}{\tau_2 - \tau_1}.
\end{eqnarray}
As R.H.S. of Eq.~\eqref{eq22neg1:Thm4} tends to zero as $\underline{\mathcal{A}}$ tends to $\underline{\hat{\mathcal{A}}}$, the continuity of Eq.~\eqref{eq3:Thm4 20} is proved. 

The last piece of this proof is to show Eq.~\eqref{eq:Thm4 21}. By setting $\mathcal{Y} \define \mathcal{Y}''$ and $\mathcal{X} \define \psi(\underline{\mathcal{A}},\mathcal{Y}')$ into the inequality provided by Eq.~\eqref{eq:10-2}, we obtain
\begin{eqnarray}\label{eq:22}
\left\Vert\psi(\underline{\mathcal{A}},\mathcal{Y}') - \psi(\underline{\mathcal{A}},\mathcal{Y}'')\right\Vert_{\rho} \leq \frac{\left\Vert \mathcal{Y}' - \mathcal{Y}''\right\Vert_{\rho}}{\tau_2 - \tau_1},
\end{eqnarray}
where $(\underline{\mathcal{A}},\mathcal{Y}'), (\underline{\mathcal{A}},\mathcal{Y}'') \in V$. Choose $V_{\epsilon} \subset V$ so that $(\underline{\mathcal{A}},\mathcal{Y})  \in V_{\epsilon}$ implies $(\underline{\mathcal{A}},\psi(\underline{\mathcal{A}},\mathcal{Y})) \in W_{\epsilon\tau_2(\tau_2 - \tau_1)}$, where $W_{\epsilon}$ is the open set required to have Eq.~\eqref{eq2:lma:Thm 3, req iii} in  Lemma~\ref{lma:Thm 3, req iii}. By selecting two elements $(\underline{\mathcal{A}},\mathcal{Y}'), (\underline{\mathcal{A}},\mathcal{Y}'') \in V_{\epsilon}$, and setting $\mathcal{X}' \define \psi(\underline{\mathcal{A}},\mathcal{Y}'), \mathcal{X}'' \define \psi(\underline{\mathcal{A}},\mathcal{Y}'')$ in Eq.~\eqref{eq2:lma:Thm 3, req iii} of Lemma~\ref{lma:Thm 3, req iii}, we have
\begin{eqnarray}
\mathcal{Y}' - \mathcal{Y}''&=& \acute{\mathcal{A}}(\psi(\underline{\mathcal{A}},\mathcal{Y}') - \psi(\underline{\mathcal{A}},\mathcal{Y}'')) + \epsilon\tau_2(\tau_2 - \tau_1)\left\Vert\psi(\underline{\mathcal{A}},\mathcal{Y}') - \psi(\underline{\mathcal{A}},\mathcal{Y}'')\right\Vert_{\rho}\mathcal{B}_{S_{\mbox{\tiny PD}}}\nonumber \\
&\subset_1&  \acute{\mathcal{A}}(\psi(\underline{\mathcal{A}},\mathcal{Y}') - \psi(\underline{\mathcal{A}},\mathcal{Y}'')) + \epsilon\tau_2\left\Vert \mathcal{Y}' - \mathcal{Y}''\right\Vert_{\rho}\mathcal{B}_{S_{\mbox{\tiny PD}}},
\end{eqnarray}
where $\mathcal{B}_{S_{\mbox{\tiny PD}}}$ is a bounded norm (by 1) tensor in the space $S_{\mbox{\tiny PD}}$ and $\subset_1$ is obtained by applying the relation provided by Eq.~\eqref{eq:22}. This means that 
we can find an tenor $\acute{\acute{\mathcal{A}}} \in \mathfrak{A}$, such that 
\begin{eqnarray}\label{eq:23-1}
\mathcal{Y}' - \mathcal{Y}''&=&\acute{\acute{\mathcal{A}}}(\psi(\underline{\mathcal{A}},\mathcal{Y}') - \psi(\underline{\mathcal{A}},\mathcal{Y}'')) + \epsilon\tau_2\left\Vert \mathcal{Y}' - \mathcal{Y}''\right\Vert_{\rho}\mathcal{B}_{S_{\mbox{\tiny PD}}}.
\end{eqnarray}
By taking the inverse of $\acute{\acute{\mathcal{A}}}$ (since all tensors in $\mathfrak{A}$ are injective) in Eq.~\eqref{eq:23-1} and using the fact that $\acute{\acute{\mathcal{A}}}^{-1}\mathbb{B}_{S_{\mbox{\tiny PD}}} \subset \frac{1}{\tau_2}\mathbb{B}_{S_{\mbox{\tiny PD}}}$, we have
\begin{eqnarray}\label{eq:23}
\psi(\underline{\mathcal{A}},\mathcal{Y}') - \psi(\underline{\mathcal{A}},\mathcal{Y}'')
= \acute{\acute{\mathcal{A}}}^{-1}\star_N (\mathcal{Y}' - \mathcal{Y}'')+ \epsilon \left\Vert \mathcal{Y}' - \mathcal{Y}'' \right\Vert_{\rho}\mathcal{B}.
\end{eqnarray}
$\hfill \Box$

\section{Tail bounds for Multivariate Random Tensor with Power Mean and Karcher Mean}\label{sec:Tail bounds for Multivariate Random Tensor Power Mean and Karcher Mean}

In this section, we will establish tail bounds for multivariate random tensors with power mean and Karcher mean. Following lemmas are required to prove the main results of this section.  

\begin{lemma}\label{lma:thm3.1}
Let $\mathfrak{M}: S^k_{\mbox{\tiny PD}} \rightarrow S_{\mbox{\tiny PD}}$ be a map of \emph{$k$-variable tensor mean} and $\sigma$ is a p.m.i. operator mean. If we have 
\begin{eqnarray}\label{eq1-1:lma:thm3.1}
\mathfrak{M}(\mathcal{A}_1^p,\cdots,\mathcal{A}_k^p) &\succeq& \lambda^{p-1}_{\min}(\mathfrak{M}(\mathcal{A}_1,\cdots,\mathcal{A}_k)) \mathfrak{M}(\mathcal{A}_1,\cdots,\mathcal{A}_k),
\end{eqnarray}
where $p \geq 1$, then, 
\begin{eqnarray}\label{eq1-2:lma:thm3.1}
\mathfrak{M}_{\sigma}(\mathcal{A}_1^p,\cdots,\mathcal{A}_k^p) &\succeq& \lambda^{p-1}_{\min}(\mathfrak{M}_{\sigma}(\mathcal{A}_1,\cdots,\mathcal{A}_k)) \mathfrak{M}_{\sigma}(\mathcal{A}_1,\cdots,\mathcal{A}_k).
\end{eqnarray}

On ther other hand, if we have 
\begin{eqnarray}\label{eq2-1:lma:thm3.1}
\mathfrak{M}(\mathcal{A}_1^p,\cdots,\mathcal{A}_k^p) &\preceq& \lambda^{p-1}_{\min}(\mathfrak{M}(\mathcal{A}_1,\cdots,\mathcal{A}_k)) \mathfrak{M}(\mathcal{A}_1,\cdots,\mathcal{A}_k),
\end{eqnarray}
where $0 < p \leq 1$, then, 
\begin{eqnarray}\label{eq2-2:lma:thm3.1}
\mathfrak{M}_{\sigma}(\mathcal{A}_1^p,\cdots,\mathcal{A}_k^p) &\preceq& \lambda^{p-1}_{\min}(\mathfrak{M}_{\sigma}(\mathcal{A}_1,\cdots,\mathcal{A}_k)) \mathfrak{M}_{\sigma}(\mathcal{A}_1,\cdots,\mathcal{A}_k).
\end{eqnarray}
\end{lemma}
\textbf{Proof:}
We will prove Eq.~\eqref{eq1-2:lma:thm3.1} first. We begin with the assumption that $1 \leq p \leq 2$, define $\mathcal{X}$ as $\mathcal{X} \define \mathfrak{M}_{\sigma}(\mathcal{A}_1,\cdots,\mathcal{A}_k)$, and set $\acute{\lambda}$ by $\acute{\lambda} \define \lambda_{\min}(\mathcal{X})$. From the congruence invaraince provided by Eq.~\eqref{eq:congruence inv}, we have
\begin{eqnarray}\label{eq3:lma:thm3.1}
\mathcal{I}=\mathfrak{M}(g_{\sigma}(\mathcal{X}^{-1/2}\mathcal{A}_1 \mathcal{X}^{-1/2}),\cdots,g_{\sigma}(\mathcal{X}^{-1/2}\mathcal{A}_k \mathcal{X}^{-1/2})).
\end{eqnarray}

From (III) in Theorem~\ref{thm:2.1}, it is enough to prove
\begin{eqnarray}
\acute{\lambda}^{p-1}\mathcal{X} &\preceq& \mathfrak{M}((\acute{\lambda}^{p-1}\mathcal{X})\sigma\mathcal{A}_1^p,\cdots,(\acute{\lambda}^{p-1}\mathcal{X})\sigma\mathcal{A}_k^p),
\end{eqnarray}
and, from congruence invariance, this is equivalent to prove
\begin{eqnarray}
\mathcal{I} &\preceq& \mathfrak{M}(g_{\sigma}(\acute{\lambda}^{1-p}\mathcal{X}^{-1/2}\mathcal{A}_1^p\mathcal{X}^{-1/2}),\cdots,g_{\sigma}(\acute{\lambda}^{1-p}\mathcal{X}^{-1/2}\mathcal{A}_k^p\mathcal{X}^{-1/2})). 
\end{eqnarray}

From Hansen-Pedersen’s inequality, for $1 \leq i \leq k$, we have 
\begin{eqnarray}
g_{\sigma}(\acute{\lambda}^{1-p}\mathcal{X}^{-1/2}\mathcal{A}_i^p\mathcal{X}^{-1/2})
\succeq g_{\sigma}((\mathcal{X}^{-1/2}\mathcal{A}_i\mathcal{X}^{-1/2})^p) 
\succeq_1 g^{p}_{\sigma}(\mathcal{X}^{-1/2}\mathcal{A}_i\mathcal{X}^{-1/2}), 
\end{eqnarray}
where $\succeq_1$ comes from p.m.i. of $\sigma$. From the assumption of Eq.~\eqref{eq1-1:lma:thm3.1}, we have
\begin{eqnarray}
\lefteqn{\mathfrak{M}(g_{\sigma}(\acute{\lambda}^{1-p}\mathcal{X}^{-1/2}\mathcal{A}_1^p\mathcal{X}^{-1/2}),\cdots,
g_{\sigma}(\acute{\lambda}^{1-p}\mathcal{X}^{-1/2}\mathcal{A}_k^p\mathcal{X}^{-1/2}))
}\nonumber \\&&
\succeq \mathfrak{M}(g^p_{\sigma}(\mathcal{X}^{-1/2}\mathcal{A}_1\mathcal{X}^{-1/2}),\cdots,
g^p_{\sigma}(\mathcal{X}^{-1/2}\mathcal{A}_k\mathcal{X}^{-1/2}))  \nonumber \\
&& \succeq \lambda_{\min}^{p-1}\left(\mathfrak{M}(g_{\sigma}(\mathcal{X}^{-1/2}\mathcal{A}_1^p\mathcal{X}^{-1/2}),\cdots,g_{\sigma}(\mathcal{X}^{-1/2}\mathcal{A}_k^p\mathcal{X}^{-1/2}))\right)  \nonumber \\
&& \times \mathfrak{M}(g_{\sigma}(\mathcal{X}^{-1/2}\mathcal{A}_1^p\mathcal{X}^{-1/2}),\cdots,g_{\sigma}(\mathcal{X}^{-1/2}\mathcal{A}_k^p\mathcal{X}^{-1/2})) \succeq \lambda_{\min}^{p-1}(\mathcal{I})\mathcal{I} = \mathcal{I}.
\end{eqnarray}
Therefore, we prove Eq.~\eqref{eq1-2:lma:thm3.1} for $1 \leq p \leq 2$. For $p \geq 2$, we willuse induction. Suppose we have Eq.~\eqref{eq1-2:lma:thm3.1} for $1 \leq p \leq 2^m$, for $2^m \leq p \leq 2^{m+1}$, we can express $p =2p'$ with $1 \leq p' \leq 2^m$. Then, 
\begin{eqnarray}
\mathfrak{M}_{\sigma}(\mathcal{A}_1^p,\cdots, \mathcal{A}_k^p) &\succeq& 
\lambda_{\min}(\mathfrak{M}_{\sigma}(\mathcal{A}_1^{p'},\cdots, \mathcal{A}_k^{p'}))
\mathfrak{M}_{\sigma}(\mathcal{A}_1^{p'},\cdots, \mathcal{A}_k^{p'}) \nonumber \\
&\succeq& \lambda_{\min}(\lambda^{p'-1}_{\min}(\mathfrak{M}_{\sigma}(\mathcal{A}_1,\cdots, \mathcal{A}_k))\mathfrak{M}_{\sigma}(\mathcal{A}_1,\cdots, \mathcal{A}_k)) \nonumber \\
& &\times \lambda^{p'-1}_{\min}(\mathfrak{M}_{\sigma}(\mathcal{A}_1,\cdots, \mathcal{A}_k))\mathfrak{M}_{\sigma}(\mathcal{A}_1,\cdots, \mathcal{A}_k) \nonumber \\
&=& \lambda^{p-1}_{\min}(\mathfrak{M}_{\sigma}(\mathcal{A}_1,\cdots, \mathcal{A}_k))\mathfrak{M}_{\sigma}(\mathcal{A}_1,\cdots, \mathcal{A}_k).
\end{eqnarray}

For the proof of Eq.~\eqref{eq2-2:lma:thm3.1}, the proof will be similar to the proof of Eq.~\eqref{eq1-2:lma:thm3.1}by reversing all inequalities and using Hansen-Pedersen’s inequality for the operator monotone function $x^p$ for $0 \leq p \leq 1$ which is  
\begin{eqnarray}
g_{\sigma}(\acute{\lambda}^{1-p}\mathcal{X}^{-1/2}\mathcal{A}_i^p\mathcal{X}^{-1/2})
\preceq g_{\sigma}((\mathcal{X}^{-1/2}\mathcal{A}_i\mathcal{X}^{-1/2})^p) 
\preceq g^{p}_{\sigma}(\mathcal{X}^{-1/2}\mathcal{A}_i\mathcal{X}^{-1/2}), 
\end{eqnarray}
where $1 \leq i \leq k$.
$\hfill \Box$

From the adjoint definition provided by Eq.~\eqref{eq:adjoint def}, if we have Eq.~\eqref{eq1-1:lma:thm3.1}, the adjoint of $\mathfrak{M}$ will satsify
\begin{eqnarray}\label{eq:3.2}
\mathfrak{M}^*(\mathcal{A}_1^p,\cdots,\mathcal{A}_k^p) &\preceq& \lambda^{p-1}_{\max}(\mathfrak{M}^*(\mathcal{A}_1,\cdots,\mathcal{A}_k)) \mathfrak{M}^*(\mathcal{A}_1,\cdots,\mathcal{A}_k),
\end{eqnarray}
where $p \geq 1$. Similarly,  if we have Eq.~\eqref{eq2-1:lma:thm3.1}, the adjoint of $\mathfrak{M}$ will satsify
\begin{eqnarray}\label{eq:3.4}
\mathfrak{M}^*(\mathcal{A}_1^p,\cdots,\mathcal{A}_k^p) &\succeq& \lambda^{p-1}_{\max}(\mathfrak{M}^*(\mathcal{A}_1,\cdots,\mathcal{A}_k)) \mathfrak{M}^*(\mathcal{A}_1,\cdots,\mathcal{A}_k),
\end{eqnarray}
where $0 < p \leq 1$. 

Next lemma is about validities of Eq.~\eqref{eq1-2:lma:thm3.1}, Eq.~\eqref{eq2-2:lma:thm3.1}, Eq.~\eqref{eq:3.2} and Eq.~\eqref{eq:3.4} for special tensor means $\mathfrak{M}$.
\begin{lemma}\label{lma:3.2}
(1) Inequalities given by Eq.~\eqref{eq1-2:lma:thm3.1} and Eq.~\eqref{eq2-2:lma:thm3.1} are valid for $\mathfrak{M}=\mathfrak{M}_{A, \bm{w}}$. \\
(2) Inequalities given by Eq.~\eqref{eq:3.2} and Eq.~\eqref{eq:3.4} are valid for $\mathfrak{M}^*=\mathfrak{M}_{H, \bm{w}}$. \\
\end{lemma}
\textbf{Proof:}
Let $\mathcal{A}_i \in S_{\mbox{\tiny PD}}$ for $1 \leq i \leq k$ and assume $1 \leq p \leq 2$, by the operator convexity of $x^p$, we have
\begin{eqnarray}\label{eq1:lma:3.2}
\sum\limits_{i=1}^k w_i \mathcal{A}_i^p &\succeq&\left(\sum\limits_{i=1}^k w_i \mathcal{A}_i\right)^p \nonumber \\
&\succeq& \lambda_{\min}^{p-1}\left(\sum\limits_{i=1}^k w_i \mathcal{A}_i\right)\sum\limits_{i=1}^k w_i \mathcal{A}_i.
\end{eqnarray} 
From the induction argument used in Lemma~\ref{lma:thm3.1}, the inequality provided by Eq.~\eqref{eq1:lma:3.2} can be extended to $p \geq 1$. 

On the other hand, given $0 < p \leq 1$, we have
\begin{eqnarray}\label{eq2:lma:3.2}
\sum\limits_{i=1}^k w_i \mathcal{A}_i^p &\preceq_1&\left(\sum\limits_{i=1}^k w_i \mathcal{A}_i\right)^p \nonumber \\
&\succeq_2& \lambda_{\min}^{p-1}\left(\sum\limits_{i=1}^k w_i \mathcal{A}_i\right)\sum\limits_{i=1}^k w_i \mathcal{A}_i,
\end{eqnarray} 
where $\preceq_1$ comes from the operator concavity of $x^p$, and $\preceq_2$ comes from the following relation
\begin{eqnarray}
\mathcal{X} \succeq \lambda^{1-p}_{\min}(\mathcal{X}) \mathcal{X}^p,~\mbox{for $\mathcal{X} \in S_{\mbox{\tiny PD}}$}.
\end{eqnarray}
From Eq.~\eqref{eq1:lma:3.2} and Eq.~\eqref{eq2:lma:3.2}, the part (1) of this Lemma is proved. The part (2) of this Lemma is also true since $\mathfrak{M}_{H,\bm{w}} = \mathfrak{M}^*_{A,\bm{w}}$.
$\hfill \Box$

We are ready to present the main theorem of this section.
\begin{theorem}\label{thm:cor:3.3}
For $p \geq 1$, $0 < q \leq 1$, $r \geq 1$, random tensors $\mathcal{X}_i \in S_{\mbox{\tiny PD}}$ for $1 \leq i \leq k$ and a deterministic PD tensor $\mathcal{C}$, we have
\begin{eqnarray}\label{eq1:thm:cor:3.3}
\mathrm{Pr}\left(\lambda_{\min}^{p-1}\left(\mathfrak{P}_{\bm{w},q}(\mathcal{X}_1,\cdots,\mathcal{X}_k)\right)\mathfrak{P}_{\bm{w},q}(\mathcal{X}_1,\cdots,\mathcal{X}_k) \npreceq \mathcal{C}\right)
\leq \mathrm{Tr}\left(\mathbb{E}\left[\left(\mathfrak{P}_{\bm{w},q}(\mathcal{X}_1^p,\cdots,\mathcal{X}_k^p)\right)^r\right]\star_N \mathcal{C}^{-1}\right);
\end{eqnarray}
Also, we have
\begin{eqnarray}\label{eq2:thm:cor:3.3}
\lefteqn{\mathrm{Pr}\left(\mathfrak{P}_{\bm{w},-q}(\mathcal{X}_1^p,\cdots,\mathcal{X}_k^p)\npreceq \mathcal{C}\right)
\leq} \nonumber \\
&&\mathrm{Tr}\left(\mathbb{E}\left[\left(\lambda^{p-1}_{\max}(\mathfrak{P}_{\bm{w},-q}(\mathcal{X}_1,\cdots,\mathcal{X}_k))\mathfrak{P}_{\bm{w},-q}(\mathcal{X}_1,\cdots,\mathcal{X}_k)\right)^r\right]\star_N \mathcal{C}^{-1}\right).
\end{eqnarray}

On the other hand, for $0 < p \leq 1$, $0 < q \leq 1$, random tensors $\mathcal{X}_i \in S_{\mbox{\tiny PD}}$ for $1 \leq i \leq k$ and a deterministic PD tensor $\mathcal{C}$, we have
\begin{eqnarray}\label{eq3:thm:cor:3.3}
\mathrm{Pr}\left(\mathfrak{P}_{\bm{w},q}(\mathcal{X}_1^p,\cdots,\mathcal{X}_k^p)\npreceq \mathcal{C}\right)
\leq \mathrm{Tr}\left(\mathbb{E}\left[\left(\lambda^{p-1}_{\min}(\mathfrak{P}_{\bm{w},q}(\mathcal{X}_1,\cdots,\mathcal{X}_k))\mathfrak{P}_{\bm{w},q}(\mathcal{X}_1,\cdots,\mathcal{X}_k)\right)^r\right]\star_N \mathcal{C}^{-1}\right);
\end{eqnarray}
Also, we have
\begin{eqnarray}\label{eq4:thm:cor:3.3}
\lefteqn{\mathrm{Pr}\left(\lambda^{p-1}_{\max}(\mathfrak{P}_{\bm{w},-q}(\mathcal{X}_1,\cdots,\mathcal{X}_k))\mathfrak{P}_{\bm{w},-q}(\mathcal{X}_1,\cdots,\mathcal{X}_k)\npreceq \mathcal{C}\right)}\nonumber \\
&\leq& \mathrm{Tr}\left(\mathbb{E}\left[\left(\mathfrak{P}_{\bm{w},-q}(\mathcal{X}_1^p,\cdots,\mathcal{X}_k^p)\right)^r\right]\star_N \mathcal{C}^{-1}\right).
\end{eqnarray}
\end{theorem}
\textbf{Proof:}
Because the operation is $\#_q$ is p.m.i., and we have $\mathfrak{P}_{\bm{w},q}(\mathcal{X}_1,\cdots,\mathcal{X}_k)$ from the solution of Eq.~\eqref{eq1:exp:weighted power mean}, and $\mathfrak{P}_{\bm{w},-q}=\mathfrak{P}^*_{\bm{w},q}$ for $0 < q \leq 1$, from Lemma~\ref{lma:thm3.1} and Lemma~\ref{lma:3.2}, we obtain the following:
\begin{eqnarray}
\mathfrak{P}_{\bm{w},q}(\mathcal{X}_1^p,\cdots,\mathcal{X}_k^p) &\succeq& 
\lambda_{\min}^{p-1}\left(\mathfrak{P}_{\bm{w},q}(\mathcal{X}_1,\cdots,\mathcal{X}_k)\right)\mathfrak{P}_{\bm{w},q}(\mathcal{X}_1,\cdots,\mathcal{X}_k), \nonumber \\
\mathfrak{P}_{\bm{w},-q}(\mathcal{X}_1^p,\cdots,\mathcal{X}_k^p) &\preceq&
\lambda^{p-1}_{\max}(\mathfrak{P}_{\bm{w},-q}(\mathcal{X}_1,\cdots,\mathcal{X}_k))\mathfrak{P}_{\bm{w},-q}(\mathcal{X}_1,\cdots,\mathcal{X}_k),
\end{eqnarray}
where $p \geq 1$.

For $0 < p \leq 1$, we have
\begin{eqnarray}
\mathfrak{P}_{\bm{w},q}(\mathcal{X}_1^p,\cdots,\mathcal{X}_k^p) &\preceq& 
\lambda_{\min}^{p-1}\left(\mathfrak{P}_{\bm{w},q}(\mathcal{X}_1,\cdots,\mathcal{X}_k)\right)\mathfrak{P}_{\bm{w},q}(\mathcal{X}_1,\cdots,\mathcal{X}_k), \nonumber \\
\mathfrak{P}_{\bm{w},-q}(\mathcal{X}_1^p,\cdots,\mathcal{X}_k^p) &\succeq&
\lambda^{p-1}_{\max}(\mathfrak{P}_{\bm{w},-q}(\mathcal{X}_1,\cdots,\mathcal{X}_k))\mathfrak{P}_{\bm{w},-q}(\mathcal{X}_1,\cdots,\mathcal{X}_k).
\end{eqnarray}

Apply Lemma 3 in~\cite{chang2023randomI} to above inequalities, we have desired results.
$\hfill \Box$

From Theorem~\ref{thm:cor:3.3}, we can have the following Corollary about Karcher mean $\mathfrak{G}_{\bm{w}}(\mathcal{A}_1,\cdots,\mathcal{A}_k)$.
\begin{corollary}\label{cor:3.4}
Given $p \geq 1$, $r \geq 1$, random tensors $\mathcal{X}_i \in S_{\mbox{\tiny PD}}$ for $1 \leq i \leq k$ and a deterministic PD tensor $\mathcal{C}$, we have
\begin{eqnarray}\label{eq1:cor:3.4}
\mathrm{Pr}\left(\lambda^{p-1}_{\min}(\mathfrak{G}_{\bm{w}}(\mathcal{A}_1,\cdots,\mathcal{A}_k))\mathfrak{G}_{\bm{w}}(\mathcal{A}_1,\cdots,\mathcal{A}_k)\npreceq \mathcal{C}\right)
\leq \mathrm{Tr}\left(\mathbb{E}\left[\left(\mathfrak{G}_{\bm{w}}(\mathcal{A}^p_1,\cdots,\mathcal{A}^p_k)\right)^r\right]\star_N \mathcal{C}^{-1}\right);
\end{eqnarray}
Also, we have
\begin{eqnarray}\label{eq2:cor:3.4}
\lefteqn{\mathrm{Pr}\left(\mathfrak{G}_{\bm{w}}(\mathcal{A}^p_1,\cdots,\mathcal{A}^p_k)\npreceq \mathcal{C}\right)
\leq} \nonumber \\
&&\mathrm{Tr}\left(\mathbb{E}\left[\left(\lambda^{p-1}_{\max}(\mathfrak{G}_{\bm{w}}(\mathcal{X}_1,\cdots,\mathcal{X}_k))\mathfrak{G}_{\bm{w}}(\mathcal{X}_1,\cdots,\mathcal{X}_k)\right)^r\right]\star_N \mathcal{C}^{-1}\right).
\end{eqnarray}

On the other hand, given $0 < p \leq 1$, $r \geq 1$, random tensors $\mathcal{X}_i \in S_{\mbox{\tiny PD}}$ for $1 \leq i \leq k$ and a deterministic PD tensor $\mathcal{C}$, we have
\begin{eqnarray}\label{eq3:cor:3.4}
\mathrm{Pr}\left(\lambda^{p-1}_{\max}(\mathfrak{G}_{\bm{w}}(\mathcal{A}_1,\cdots,\mathcal{A}_k))\mathfrak{G}_{\bm{w}}(\mathcal{A}_1,\cdots,\mathcal{A}_k)\npreceq \mathcal{C}\right)
\leq \mathrm{Tr}\left(\mathbb{E}\left[\left(\mathfrak{G}_{\bm{w}}(\mathcal{A}^p_1,\cdots,\mathcal{A}^p_k)\right)^r\right]\star_N \mathcal{C}^{-1}\right);
\end{eqnarray}
Also, we have
\begin{eqnarray}\label{eq4:cor:3.4}
\lefteqn{\mathrm{Pr}\left(\mathfrak{G}_{\bm{w}}(\mathcal{A}^p_1,\cdots,\mathcal{A}^p_k)\npreceq \mathcal{C}\right)
\leq} \nonumber \\
&&\mathrm{Tr}\left(\mathbb{E}\left[\left(\lambda^{p-1}_{\min}(\mathfrak{G}_{\bm{w}}(\mathcal{X}_1,\cdots,\mathcal{X}_k))\mathfrak{G}_{\bm{w}}(\mathcal{X}_1,\cdots,\mathcal{X}_k)\right)^r\right]\star_N \mathcal{C}^{-1}\right).
\end{eqnarray}
\end{corollary}
\textbf{Proof:}
As $q \searrow 0$, we have
\begin{eqnarray}\label{eq5:cor:3.4}
\lambda_{\min}(\mathfrak{P}_{\bm{w},q}(\mathcal{X}_1,\cdots,\mathcal{X}_k)) \searrow 
\lambda_{\min}(\mathfrak{G}_{\bm{w}}(\mathcal{X}_1,\cdots,\mathcal{X}_k)),
\end{eqnarray}
and
\begin{eqnarray}\label{eq6:cor:3.4}
\lambda_{\max}(\mathfrak{P}_{\bm{w},-q}(\mathcal{X}_1,\cdots,\mathcal{X}_k)) \nearrow 
\lambda_{\max}(\mathfrak{G}_{\bm{w}}(\mathcal{X}_1,\cdots,\mathcal{X}_k)).
\end{eqnarray}
By taking the limits of those inequalities in the proof in Theorem~\ref{thm:cor:3.3}, we have
\begin{eqnarray}\label{eq7:cor:3.4}
\lambda^{p-1}_{\min}(\mathfrak{G}_{\bm{w}}(\mathcal{A}_1,\cdots,\mathcal{A}_k))\mathfrak{G}_{\bm{w}}(\mathcal{A}_1,\cdots,\mathcal{A}_k) &\preceq&\mathfrak{G}_{\bm{w}}(\mathcal{A}^p_1,\cdots,\mathcal{A}^p_k) \nonumber \\
&\preceq& \lambda^{p-1}_{\max}(\mathfrak{G}_{\bm{w}}(\mathcal{X}_1,\cdots,\mathcal{X}_k))\mathfrak{G}_{\bm{w}}(\mathcal{X}_1,\cdots,\mathcal{X}_k),
\end{eqnarray}
where $p \geq 1$; and
\begin{eqnarray}\label{eq8:cor:3.4}
\lambda^{p-1}_{\max}(\mathfrak{G}_{\bm{w}}(\mathcal{A}_1,\cdots,\mathcal{A}_k))\mathfrak{G}_{\bm{w}}(\mathcal{A}_1,\cdots,\mathcal{A}_k) &\preceq& \mathfrak{G}_{\bm{w}}(\mathcal{A}^p_1,\cdots,\mathcal{A}^p_k) \nonumber \\
&\preceq& \lambda^{p-1}_{\min}(\mathfrak{G}_{\bm{w}}(\mathcal{X}_1,\cdots,\mathcal{X}_k))\mathfrak{G}_{\bm{w}}(\mathcal{X}_1,\cdots,\mathcal{X}_k),
\end{eqnarray}
where $0 < p \leq 1$.

Apply Lemma 3 in~\cite{chang2023randomI} to above inequalities provided by Eq.~\eqref{eq7:cor:3.4} and~\eqref{eq8:cor:3.4}, we have desired results.
$\hfill \Box$

\section{Tail bounds for Multivariate Random Tensor Deformed Means Based on Ando-Hiai Type Inequalities}\label{sec:Tail bounds for Multivariate Random Tensor Deformed Means Based on Ando-Hiai Type Inequalities}

In this section, we will generalize tail bounds for multivariate random tensor for the deformed means $\mathfrak{M}_{\sigma_{1/p}}$ and $\mathfrak{M}_{\sigma_{p}}$, where  $\sigma_{1/p}$ and $\sigma_{p}$ are operator means with respect to the representing functions $g_{\sigma}(x^{1/p})$ and $g_{\sigma}(x^p)$, respectively. We begin this section with the following two lemmas. The first lemma is about $\mathfrak{M}_{\sigma_{1/p}}$. 

\begin{lemma}\label{lma:thm:4.1}
Let $\mathfrak{M}: S^k_{\mbox{\tiny PD}} \rightarrow S_{\mbox{\tiny PD}}$ be a map of \emph{$k$-variable tensor mean} and $\sigma_{1/p}$ be the operator mean with the representing function $g_{\sigma}(x^{1/p})$ for $p \geq 1$. Then 
\begin{eqnarray}\label{eq1:lma:thm:4.1}
\lambda^{p-1}_{\min}(\mathfrak{M}_{\sigma}(\mathcal{A}_1,\cdots,\mathcal{A}_k))
\mathfrak{M}_{\sigma}(\mathcal{A}_1,\cdots,\mathcal{A}_k) &\preceq_1& 
\mathfrak{M}_{\sigma_{1/p}}(\mathcal{A}^p_1,\cdots,\mathcal{A}^p_k) \nonumber \\
&\preceq_2& 
\lambda^{p-1}_{\max}(\mathfrak{M}_{\sigma}(\mathcal{A}_1,\cdots,\mathcal{A}_k))
\mathfrak{M}_{\sigma}(\mathcal{A}_1,\cdots,\mathcal{A}_k).
\end{eqnarray}
\end{lemma}
\textbf{Proof:}
We will prove the inequality $\preceq_2$ in Eq.~\eqref{eq1:lma:thm:4.1} first. Let $\mathcal{X} \define \mathfrak{M}_{\sigma}(\mathcal{A}_1,\cdots,\mathcal{A}_k)$, and the statement (IV) in Theorem~\ref{thm:2.1}, it is enough for us to prove the following 
\begin{eqnarray}\label{eq2:lma:thm:4.1}
\lambda^{p-1}_{\max}(\mathcal{X})\mathcal{X}&\succeq&\mathfrak{M}((\lambda^{p-1}_{\max}(\mathcal{X})\mathcal{X})\sigma_{1/p}\mathcal{A}^p_1,\cdots,(\lambda^{p-1}_{\max}(\mathcal{X})\mathcal{X})\sigma_{1/p}\mathcal{A}^p_k),
\end{eqnarray}
which can also be expressed by
\begin{eqnarray}\label{eq3:lma:thm:4.1}
\mathcal{I}&\succeq&\mathfrak{M}(g_{\sigma}(\lambda^{\frac{1}{p}-1}_{\max}(\mathcal{X})(\mathcal{X}^{-1/2} \mathcal{A}^p_1\mathcal{X}^{-1/2})^{1/p}),\cdots,g_{\sigma}(\lambda^{\frac{1}{p}-1}_{\max}(\mathcal{X})(\mathcal{X}^{-1/2} \mathcal{A}^p_k\mathcal{X}^{-1/2})^{1/p})).
\end{eqnarray}

From Hansen-Pedersen’s inequality, for $1 \leq i \leq k$, we have 
\begin{eqnarray}\label{eq4:lma:thm:4.1}
\lambda^{\frac{1}{p}-1}_{\max}(\mathcal{X}^{-1/2} \mathcal{A}^p_i \mathcal{X}^{-1/2})^{1/p}
&\preceq&\mathcal{X}^{-1/2} \mathcal{A}_i \mathcal{X}^{-1/2},
\end{eqnarray}
then, from the monotonicity of $g_{\sigma}$ and the congruence property of $\mathfrak{M}$, we have
\begin{eqnarray}\label{eq5:lma:thm:4.1}
\lefteqn{\mathfrak{M}(g_{\sigma}(\lambda^{\frac{1}{p}-1}_{\max}(\mathcal{X})(\mathcal{X}^{-1/2} \mathcal{A}^p_1\mathcal{X}^{-1/2})^{1/p}),\cdots,g_{\sigma}(\lambda^{\frac{1}{p}-1}_{\max}(\mathcal{X})(\mathcal{X}^{-1/2} \mathcal{A}^p_k\mathcal{X}^{-1/2})^{1/p}))}&& \nonumber \\
&\preceq& \mathfrak{M}(g_{\sigma}((\mathcal{X}^{-1/2} \mathcal{A}^p_1\mathcal{X}^{-1/2})^{1/p}),\cdots,g_{\sigma}((\mathcal{X}^{-1/2} \mathcal{A}^p_k\mathcal{X}^{-1/2})^{1/p})) = \mathcal{I}.
\end{eqnarray}
This establishes the inequality $\preceq_2$ of Eq.~\eqref{eq1:lma:thm:4.1}.

The inequality $\preceq_1$ of Eq.~\eqref{eq1:lma:thm:4.1} can be derived from $\preceq_2$ in Eq.~\eqref{eq1:lma:thm:4.1}. By replacing $\mathfrak{M}, \sigma$ and $\mathcal{A}_i$ for $1 \leq i \leq k$ in those terms in both sides of $\preceq_2$ and applying the following relations: $(\mathfrak{M}^*)_{(\sigma^*)_{1/p}} = (\mathfrak{M}_{\sigma_{1/p}})^*$ and $(\mathfrak{M}^*)_{\sigma^*} = (\mathfrak{M}_{\sigma})^*$, we have
\begin{eqnarray}
\mathfrak{M}^{-1}_{\sigma_{1/p}}(\mathcal{A}^p_1,\cdots,\mathcal{A}^p_k) 
&\preceq& 
\lambda^{p-1}_{\max}(\mathfrak{M}^{-1}_{\sigma}(\mathcal{A}_1,\cdots,\mathcal{A}_k))
\mathfrak{M}^{-1}_{\sigma}(\mathcal{A}_1,\cdots,\mathcal{A}_k),
\end{eqnarray}
which is equivalent to
\begin{eqnarray}
\lambda^{p-1}_{\min}(\mathfrak{M}_{\sigma}(\mathcal{A}_1,\cdots,\mathcal{A}_k))
\mathfrak{M}_{\sigma}(\mathcal{A}_1,\cdots,\mathcal{A}_k) &\preceq& 
\mathfrak{M}_{\sigma_{1/p}}(\mathcal{A}^p_1,\cdots,\mathcal{A}^p_k) 
\end{eqnarray}
Therefore, this lemma is proved. 
$\hfill \Box$

The next lemma is about $\mathfrak{M}_{\sigma_{p}}$. 
\begin{lemma}\label{lma:thm:4.2}
Let $\mathfrak{M}: S^k_{\mbox{\tiny PD}} \rightarrow S_{\mbox{\tiny PD}}$ be a map of \emph{$k$-variable tensor mean} and $\sigma_{p}$ be the operator mean with the representing function $g_{\sigma}(x^{p})$ for $0 < p \leq 1$. Then 
\begin{eqnarray}\label{eq1:lma:thm:4.2}
\lambda^{p-1}_{\max}(\mathfrak{M}_{\sigma_p}(\mathcal{A}_1,\cdots,\mathcal{A}_k))
\mathfrak{M}_{\sigma_p}(\mathcal{A}_1,\cdots,\mathcal{A}_k) &\preceq_1& 
\mathfrak{M}_{\sigma}(\mathcal{A}^p_1,\cdots,\mathcal{A}^p_k) \nonumber \\
&\preceq_2& 
\lambda^{p-1}_{\min}(\mathfrak{M}_{\sigma_p}(\mathcal{A}_1,\cdots,\mathcal{A}_k))
\mathfrak{M}_{\sigma_p}(\mathcal{A}_1,\cdots,\mathcal{A}_k).
\end{eqnarray}
\end{lemma}
\textbf{Proof:}
We will prove the inequality $\preceq_2$ in Eq.~\eqref{eq1:lma:thm:4.2} first. Let $\mathcal{X} \define \mathfrak{M}_{\sigma}(\mathcal{A}_1,\cdots,\mathcal{A}_k)$, and the statement (IV) in Theorem~\ref{thm:2.1}, it is enough for us to prove the following 
\begin{eqnarray}\label{eq2:lma:thm:4.2}
\lambda^{p-1}_{\min}(\mathcal{X})\mathcal{X}&\succeq&\mathfrak{M}((\lambda^{p-1}_{\min}(\mathcal{X})\mathcal{X})\sigma\mathcal{A}^p_1,\cdots,(\lambda^{p-1}_{\min}(\mathcal{X})\mathcal{X})\sigma\mathcal{A}^p_k), 
\end{eqnarray}
which can also be expressed by
\begin{eqnarray}\label{eq3:lma:thm:4.2}
\mathcal{I}&\succeq&\mathfrak{M}(g_{\sigma}(\lambda^{1-p}_{\min}(\mathcal{X})\mathcal{X}^{-1/2} \mathcal{A}^p_1\mathcal{X}^{-1/2}),\cdots,g_{\sigma}(\lambda^{1-p}_{\min}(\mathcal{X})\mathcal{X}^{-1/2} \mathcal{A}^p_k\mathcal{X}^{-1/2})^{1/p}).
\end{eqnarray}

From Hansen-Pedersen’s inequality, for $1 \leq i \leq k$, we have 
\begin{eqnarray}\label{eq4:lma:thm:4.2}
\lambda^{1-p}_{\min}(\mathcal{X})\mathcal{X}^{-1/2} \mathcal{A}^p_i \mathcal{X}^{-1/2}
&\preceq&(\mathcal{X}^{-1/2} \mathcal{A}_i \mathcal{X}^{-1/2})^p,
\end{eqnarray}
then, from the monotonicity of $g_{\sigma}$ and the congruence property of $\mathfrak{M}$, we have
\begin{eqnarray}\label{eq5:lma:thm:4.2}
\lefteqn{\mathfrak{M}(g_{\sigma}(\lambda^{1-p}_{\min}(\mathcal{X})\mathcal{X}^{-1/2} \mathcal{A}^p_1 \mathcal{X}^{-1/2}),\cdots,g_{\sigma}(\lambda^{1-p}_{\min}(\mathcal{X})\mathcal{X}^{-1/2} \mathcal{A}^p_k \mathcal{X}^{-1/2}))}&& \nonumber \\
&\preceq& \mathfrak{M}(g_{\sigma}((\mathcal{X}^{-1/2} \mathcal{A}_1\mathcal{X}^{-1/2})^{p}),\cdots,g_{\sigma}((\mathcal{X}^{-1/2} \mathcal{A}_k\mathcal{X}^{-1/2})^{p})) \nonumber \\
&=& \mathfrak{M}(g_{\sigma_p}(\mathcal{X}^{-1/2} \mathcal{A}_1\mathcal{X}^{-1/2}),\cdots,g_{\sigma_p}(\mathcal{X}^{-1/2} \mathcal{A}_k\mathcal{X}^{-1/2})) = \mathcal{I}.
\end{eqnarray}
Therefore, we have proved the inequality $\preceq_2$ in Eq.~\eqref{eq1:lma:thm:4.2}. 

The inequality $\preceq_1$ of Eq.~\eqref{eq1:lma:thm:4.2} can be derived from $\preceq_2$ in Eq.~\eqref{eq1:lma:thm:4.2}. By replacing $\mathfrak{M}, \sigma$ and $\mathcal{A}_i$ for $1 \leq i \leq k$ in those terms in both sides of $\preceq_2$ and applying the following relations: $(\mathfrak{M}^*)_{(\sigma^*)_{p}} = (\mathfrak{M}_{\sigma_{p}})^*$ and $(\mathfrak{M}^*)_{\sigma^*} = (\mathfrak{M}_{\sigma})^*$, we have
\begin{eqnarray}
\mathfrak{M}^{-1}_{\sigma}(\mathcal{A}^p_1,\cdots,\mathcal{A}^p_k) 
&\preceq& 
\lambda^{p-1}_{\min}(\mathfrak{M}^{-1}_{\sigma_p}(\mathcal{A}_1,\cdots,\mathcal{A}_k))
\mathfrak{M}^{-1}_{\sigma_p}(\mathcal{A}_1,\cdots,\mathcal{A}_k),
\end{eqnarray}
which is equivalent to
\begin{eqnarray}
\lambda^{p-1}_{\max}(\mathfrak{M}_{\sigma_p}(\mathcal{A}_1,\cdots,\mathcal{A}_k))
\mathfrak{M}_{\sigma_p}(\mathcal{A}_1,\cdots,\mathcal{A}_k) &\preceq& 
\mathfrak{M}_{\sigma}(\mathcal{A}^p_1,\cdots,\mathcal{A}^p_k).
\end{eqnarray}
Therefore, this lemma is proved. 
$\hfill \Box$

We are ready to present the main theorem of this section.
\begin{theorem}\label{thm:cor:4.4}
For $p \geq 1$, $-1 \leq q \leq 1$, $r \geq 1$, random tensors $\mathcal{X}_i \in S_{\mbox{\tiny PD}}$ for $1 \leq i \leq k$ and a deterministic PD tensor $\mathcal{C}$, we have
\begin{eqnarray}\label{eq1:thm:cor:4.4}
\lefteqn{\mathrm{Pr}\left(\lambda_{\min}^{p-1}\left(\mathfrak{P}_{\bm{w},q}(\mathcal{X}_1,\cdots,\mathcal{X}_k)\right)\mathfrak{P}_{\bm{w},q}(\mathcal{X}_1,\cdots,\mathcal{X}_k) \npreceq \mathcal{C}\right)
\leq} \nonumber \\
&& \mathrm{Tr}\left(\mathbb{E}\left[\left(\mathfrak{P}_{\bm{w},q/p}(\mathcal{X}_1^p,\cdots,\mathcal{X}_k^p)\right)^r\right]\star_N \mathcal{C}^{-1}\right);
\end{eqnarray}
Also, we have
\begin{eqnarray}\label{eq2:thm:cor:4.4}
\lefteqn{\mathrm{Pr}\left(\mathfrak{P}_{\bm{w},q/p}(\mathcal{X}_1^p,\cdots,\mathcal{X}_k^p)\npreceq \mathcal{C}\right)
\leq} \nonumber \\
&&\mathrm{Tr}\left(\mathbb{E}\left[\left(\lambda^{p-1}_{\max}(\mathfrak{P}_{\bm{w},q}(\mathcal{X}_1,\cdots,\mathcal{X}_k))\mathfrak{P}_{\bm{w},q}(\mathcal{X}_1,\cdots,\mathcal{X}_k)\right)^r\right]\star_N \mathcal{C}^{-1}\right).
\end{eqnarray}

On the other hand, for $0 < p \leq 1$, $-1 \leq q \leq 1$, random tensors $\mathcal{X}_i \in S_{\mbox{\tiny PD}}$ for $1 \leq i \leq k$ and a deterministic PD tensor $\mathcal{C}$, we have
\begin{eqnarray}\label{eq3:thm:cor:4.4}
\lefteqn{\mathrm{Pr}\left(\lambda^{p-1}_{\max}(\mathfrak{P}_{\bm{w},qp}(\mathcal{X}_1,\cdots,\mathcal{X}_k))\mathfrak{P}_{\bm{w},qp}(\mathcal{X}_1,\cdots,\mathcal{X}_k)\npreceq \mathcal{C}\right)}\nonumber \\
&\leq& \mathrm{Tr}\left(\mathbb{E}\left[\left((\mathfrak{P}_{\bm{w},q}(\mathcal{X}^p_1,\cdots,\mathcal{X}^p_k))\mathfrak{P}_{\bm{w},q}(\mathcal{X}_1,\cdots,\mathcal{X}_k)\right)^r\right]\star_N \mathcal{C}^{-1}\right);
\end{eqnarray}
Also, we have
\begin{eqnarray}\label{eq4:thm:cor:4.4}
\lefteqn{\mathrm{Pr}\left(\mathfrak{P}_{\bm{w},q}(\mathcal{X}^p_1,\cdots,\mathcal{X}^p_k)\npreceq \mathcal{C}\right)}\nonumber \\
&\leq& \mathrm{Tr}\left(\mathbb{E}\left[\left(\lambda^{p-1}_{\min}(\mathfrak{P}_{\bm{w},qp}(\mathcal{X}_1,\cdots,\mathcal{X}_k))\mathfrak{P}_{\bm{w},qp}(\mathcal{X}_1,\cdots,\mathcal{X}_k)\right)^r\right]\star_N \mathcal{C}^{-1}\right).
\end{eqnarray}
\end{theorem}
\textbf{Proof:}
From Lemma~\ref{lma:thm:4.1}, we have
\begin{eqnarray}\label{eq7:cor:4.4}
\lambda^{p-1}_{\min}(\mathfrak{P}_{\bm{w},q}(\mathcal{X}_1,\cdots,\mathcal{X}_k))\mathfrak{P}_{\bm{w},q}(\mathcal{X}_1,\cdots,\mathcal{X}_k) &\preceq&\mathfrak{P}_{\bm{w},q/p}(\mathcal{X}^p_1,\cdots,\mathcal{X}^p_k) \nonumber \\
&\preceq& \lambda^{p-1}_{\max}(\mathfrak{P}_{\bm{w},q}(\mathcal{X}_1,\cdots,\mathcal{X}_k))\mathfrak{P}_{\bm{w},q}(\mathcal{X}_1,\cdots,\mathcal{X}_k),
\end{eqnarray}
where $p \geq 1$; and, from Lemma~\ref{lma:thm:4.2}, we also have
\begin{eqnarray}\label{eq8:cor:4.4}
\lambda^{p-1}_{\max}(\mathfrak{P}_{\bm{w},qp}(\mathcal{X}_1,\cdots,\mathcal{X}_k))\mathfrak{P}_{\bm{w},qp}(\mathcal{X}_1,\cdots,\mathcal{X}_k) &\preceq& \mathfrak{P}_{\bm{w},q}(\mathcal{X}^p_1,\cdots,\mathcal{X}^p_k) \nonumber \\
&\preceq& \lambda^{p-1}_{\min}(\mathfrak{P}_{\bm{w},qp}(\mathcal{X}_1,\cdots,\mathcal{X}_k))\mathfrak{P}_{\bm{w},qp}(\mathcal{X}_1,\cdots,\mathcal{X}_k),
\end{eqnarray}
where $0 < p \leq 1$.

Apply Lemma 3 in~\cite{chang2023randomI} to above inequalities provided by Eq.~\eqref{eq7:cor:4.4} and~\eqref{eq8:cor:4.4}, we have desired results.
$\hfill \Box$

\section{Tail bounds for Multivariate Random Tensor Means Based on Reverse Ando-Hiai Type Inequalities}\label{sec:Tail bounds for Multivariate Random Tensor Means Based on Reverse Ando-Hiai Type Inequalities}

Different from previous sections, we want to obtain the reverse of Ando-Hiai type inequalities via Kantorovich constants. The Kantorovich constants are defined by
\begin{eqnarray}\label{eq:Kantorovich cons def0}
K(M,m,f,p) &=& \frac{mf(M) - Mf(m)}{(p-1)(M-m)}\left(\frac{(p-1)(f(M)-f(m))}{p(mF(M)-Mf(m))}\right)^p,
\end{eqnarray}
where $p$ is a real number and $f(t)$ is a real valued continuous function defined on an real interval $[m,M]$~\cite{josip2013MondPecaricI}. If $f(t)=t^q$, the Kantorovich constants become 
\begin{eqnarray}\label{eq:Kantorovich cons def1}
K(M,m,t^q,p) &=& \frac{mM^q - Mm^q)}{(p-1)(M-m)}\left(\frac{(p-1)(M^q-m^q)}{p(mM^q-Mm^q)}\right)^p.
\end{eqnarray}
Further, if we have $p=q$, we can express the Kantorovich constants as 
\begin{eqnarray}\label{eq:Kantorovich cons def2}
K(M,m,t^p,p) &=& \frac{mM^p - Mm^p)}{(p-1)(M-m)}\left(\frac{(p-1)(M^p-m^p)}{p(mM^p-Mm^p)}\right)^p \nonumber \\
&\define& K(M,m,p). 
\end{eqnarray}
Note that $K(M,m,p) > 1$.

We have the following lemmas based on the Kantorovich constants. 
\begin{lemma}\label{lma:5.1}
Let $\mathcal{A}, \mathcal{B} \in \mathbb{C}^{I_1 \times \dots \times I_N \times I_1 \times \dots \times I_N}$ be two PD tensors satisfying $m\mathcal{I} \preceq \mathcal{A} 
\preceq M\mathcal{I}$, and $\mathcal{B}^2 \preceq \mathcal{I}$, where $0 < m < M$. Given $p > 1$, we have
\begin{eqnarray}\label{eq1:lma:5.1}
\mathcal{B}\star_N \mathcal{A}^p \star_N \mathcal{B}&\preceq& K(M,m,p)(\mathcal{B}\star_N \mathcal{A} \star_N \mathcal{B})^p,
\end{eqnarray}
where $K(M,m,p)$ is defined by Eq~\eqref{eq:Kantorovich cons def2}.
\end{lemma}
\textbf{Proof:}
From Corollary 2.12 in Book~\cite{josip2013MondPecaricI}, we have
\begin{eqnarray}\label{eq2:lma:5.1}
\sum\limits_{j=1}^k \Phi_j(\mathcal{A}_j^p)&\preceq& K(M,m,p) \left(\sum\limits_{j=1}^k \omega_j\Phi_j(\mathcal{A}_j)\right)^p,
\end{eqnarray}
where $\sum\limits_{j=1}^k \omega_j =1$. By setting $k=1$, $\mathcal{A}_j = \mathcal{A}$, $\omega_1=1$, and $\Phi(\mathcal{X}) = \mathcal{B}\star_N \mathcal{X} \star_N \mathcal{B}$, where $\mathcal{X}$ is any PD tensor, into Eq.~\eqref{eq2:lma:5.1}, we have
\begin{eqnarray}\label{eq3:lma:5.1}
\mathcal{B}\star_N \mathcal{A}^p \star_N \mathcal{B}&\preceq& K(M,m,p)(\mathcal{B}\star_N \mathcal{A} \star_N \mathcal{B})^p.
\end{eqnarray}
Therefore, we obtain the desired inequality.
$\hfill \Box$

\begin{lemma}\label{lma:eq 5.3}
Let $\mathcal{A}_i \in \mathbb{C}^{I_1 \times \dots \times I_N \times I_1 \times \dots \times I_N}$ be PD tensors satisfying $m\mathcal{I} \preceq \mathcal{A}_i \preceq M\mathcal{I}$. Given $p > 1$, we have
\begin{eqnarray}\label{eq1:lma:eq 5.3}
\sum\limits_{i=1}^k \omega_i \mathcal{A}_i^p &\preceq& K(M,m,p)\left(\sum\limits_{i=1}^k \omega_i \mathcal{A}_i\right)^p
\end{eqnarray}
where $\omega_i \geq 0$ and $\sum\limits_{i=1}^k \omega_i =1$, i.e. $(\omega_1,\cdots,\omega_k)$ is a probability vector.
\end{lemma}
\textbf{Proof:}
From Theorem 3.18 in Book~\cite{josip2013MondPecaricI}, we have
\begin{eqnarray}\label{eq2:lma:eq 5.3}
\Phi(\mathcal{A}^p)&\preceq& K(M,m,p) \left(\Phi(\mathcal{A})\right)^p,
\end{eqnarray}
where $\Phi$ is a normalized positive linear map from a bounded linear operator space over the Hilbert space $\mathrm{H}_1$, denoted by $\mathfrak{B}(\mathrm{H}_1)$ to another bounded linear operator space over the Hilbert space $\mathrm{H}_2$, denoted by $\mathfrak{B}(\mathrm{H}_2)$, $p$ is greater or equal than $1$, and $\mathcal{A}$ is a PD tensor. Applying Eq.~\eqref{eq2:lma:eq 5.3} to the normalized positive linear map $\Psi: \mathfrak{B}(\oplus_1^{k^2}\mathfrak{H}) \rightarrow \mathfrak{B}(\mathfrak{H})$ defined by $\Psi([\mathcal{A}_{i,j}]_{i,j=1}^k) = \sum\limits_{i=1}^k \omega_i A_{i,i}$, and assuming that $m \mathcal{I} \preceq A_{i,j} \preceq M \mathcal{I}$ for all $1 \leq i, j \leq k$, we have
\begin{eqnarray}\label{eq3:lma:eq 5.3}
\sum\limits_{i=1}^k \omega_i \mathcal{A}_{i,i}^p &\preceq& K(M,m,p) \left(\sum\limits_{i=1}^k \omega_i \mathcal{A}_{i,i}\right)^p.
\end{eqnarray}
We prove this Lemma by replacing $\mathcal{A}_{i,i}$ with $\mathcal{A}_{i}$. 
$\hfill \Box$

Following Theorem is about $k$-variable power means $\mathfrak{P}_{\bm{w},q}$ which can be treated as the reverse counterparts of Theorem~\ref{thm:cor:3.3}.
\begin{theorem}\label{thm:5.2}
Let $\mathcal{A}_i \in \mathbb{C}^{I_1 \times \dots \times I_N \times I_1 \times \dots \times I_N}$ be random PD tensors satisfying $m\mathcal{I} \preceq \mathcal{A}_i \preceq M\mathcal{I}$ almost surely, a deterministic PD tensor $\mathcal{C} \in \mathbb{C}^{I_1 \times \dots \times I_N \times I_1 \times \dots \times I_N}$ and $\mathcal{X} \define \mathfrak{P}_{\bm{w},q}(\mathcal{A}_1,\cdots,\mathcal{A}_k)$. We also set $\lambda_{\min}(\mathcal{X}) \define \lambda_{\min}$, and $\lambda_{\max}(\mathcal{X}) \define \lambda_{\max}$. Then, for any $q \in (0,1]$, $r \geq 1$ and $p \geq 1$, we have
\begin{eqnarray}\label{eq1:thm:5.2}
\mathrm{Pr}\left(\mathfrak{P}_{\bm{w},q}(\mathcal{A}^p_1,\cdots,\mathcal{A}^p_k) \npreceq \mathcal{C}\right)
\leq \mathrm{Tr}\left(\mathbb{E}\left[\left(\lambda_{\min}^{p-1}K_1 K_2^{1/q}
\mathfrak{P}_{\bm{w},q}(\mathcal{A}_1,\cdots,\mathcal{A}_k)\right)^r\right] \star_N \mathcal{C}^{-1}\right),
\end{eqnarray}
where $K_1 \define  K\left(M,m,p\right)$ and $K_2 \define K\left(\frac{M^q}{\lambda^q_{\min}},\frac{m^q}{\lambda^q_{\max}},p\right)$. 

On the other hand, for any $q \in [-1,0)$,  $r \geq 1$ and $p \geq 1$, we have
\begin{eqnarray}\label{eq2:thm:5.2}
\mathrm{Pr}\left(\lambda_{\max}^{p-1}K^{-1}_1 K'^{1/q}_2\mathfrak{P}_{\bm{w},q}(\mathcal{A}_1,\cdots,\mathcal{A}_k) \npreceq \mathcal{C}\right)
\leq \mathrm{Tr}\left(\mathbb{E}\left[\left(
\mathfrak{P}_{\bm{w},q}(\mathcal{A}^p_1,\cdots,\mathcal{A}^p_k)\right)^r\right] \star_N \mathcal{C}^{-1}\right),
\end{eqnarray}
where $K'_2 \define K\left(\frac{\lambda^q_{\max}}{m^q},\frac{\lambda^q_{\min}}{M^q},p\right)$. 
\end{theorem}
\textbf{Proof:}
Since $\mathfrak{P}_{\bm{w},q}^\mathrm{H} = \mathfrak{P}_{\bm{w},-q}$, $M^{-1} \mathcal{I} \preceq \mathcal{A}^{-1}_i \preceq m^{-1}\mathcal{I}$, we will have Eq.~\eqref{eq2:thm:5.2} from Eq.~\eqref{eq1:thm:5.2} by replacing $\mathcal{A}_i$ with $\mathcal{A}^{-1}_i$ and $q$ with $-q$. Therefore, it is enough to prove  Eq.~\eqref{eq1:thm:5.2} only. 

Because $\mathcal{X}=\sum\limits_{i=1}^k\omega_i \mathcal{X} \#_{q} \mathcal{A}_i$, we have
\begin{eqnarray}\label{eq3:thm:5.2}
\mathcal{I} &=& \sum\limits_{i=1}^k \omega_i \left(\mathcal{X}^{-1/2}\star_N\mathcal{A}_i\star_N\mathcal{X}^{-1/2}\right)^q;
\end{eqnarray}
and the fact that $\lambda_{\min}\mathcal{X}^{-1} \preceq \mathcal{I}$, we have 
\begin{eqnarray}\label{eq4:thm:5.2}
\lefteqn{(\lambda^{1/2}_{\min}\mathcal{X}^{-1/2}) \star_N \mathcal{A}_i^p \star_N (\lambda^{1/2}_{\min}\mathcal{X}^{-1/2}) \preceq} \nonumber \\
&&K(M,m,p)\left[(\lambda^{1/2}_{\min}\mathcal{X}^{-1/2}) \star_N \mathcal{A}_i \star_N (\lambda^{1/2}_{\min}\mathcal{X}^{-1/2})\right]^p,
\end{eqnarray}
by Lemma~\ref{lma:5.1}. Since $\lambda_{\min}(\mathcal{X}) \define \lambda_{\min}$, and $ K\left(M,m,p\right) \define K_1$, we have 
\begin{eqnarray}\label{eq5:thm:5.2}
\sum\limits_{i=1}^k \omega_i \left[\lambda_{\min}^{1-p}K^{-1}_1(\mathcal{X}^{-1/2}\mathcal{A}^{p}_{i}\mathcal{X}^{-1/2})\right]^q \preceq \sum\limits_{i=1}^k \omega_i (\mathcal{X}^{-1/2}\mathcal{A}_{i}\mathcal{X}^{-1/2})^{qp}.
\end{eqnarray}

Besides, due to $\frac{m}{\lambda_{\max}} \mathcal{I} \preceq \mathcal{X}^{-1/2} \mathcal{A}_i \mathcal{X}^{-1/2} \preceq \frac{M}{\lambda_{\min}}\mathcal{I}$ and $0 < q  \leq 1$, we also have 
\begin{eqnarray}\label{eq6:thm:5.2}
\left(\frac{m}{\lambda_{\max}}\right)^q \mathcal{I} \preceq \left(\mathcal{X}^{-1/2} \mathcal{A}_i \mathcal{X}^{-1/2}\right)^q \preceq \left(\frac{M}{\lambda_{\min}}\right)^q \mathcal{I}.
\end{eqnarray}
From Lemma~\ref{lma:eq 5.3} and Eq.~\eqref{eq5:thm:5.2}, we have
\begin{eqnarray}\label{eq7:thm:5.2}
\lefteqn{\sum\limits_{i=1}^k \omega_i \left[\lambda_{\min}^{1-p}K^{-1}_1(\mathcal{X}^{-1/2}\mathcal{A}^{p}_{i}\mathcal{X}^{-1/2})\right]^q \preceq} \nonumber \\
&& K\left(\frac{M^q}{\lambda^q_{\min}},\frac{m^q}{\lambda^q_{\max}},p\right)\left[\sum\limits_{i=1}^k \omega_i\left(\mathcal{X}^{-1/2}\mathcal{A}_{i}\mathcal{X}^{-1/2}\right)^q\right]^p.
\end{eqnarray}
By setting $K\left(\frac{M^q}{\lambda^q_{\min}},\frac{m^q}{\lambda^q_{\max}},p\right) \define K_2$ and Eq.~\eqref{eq7:thm:5.2}, we have 
\begin{eqnarray}\label{eq8:thm:5.2}
\sum\limits_{i=1}^k \omega_i \left[\left(\lambda_{\min}^{p-1}K_1 K_2^{1/q}\mathcal{X}\right)\#_{q}\mathcal{A}_i^p\right] \preceq \lambda_{\min}^{p-1} K_1 K_2^{1/q} \mathcal{X}.
\end{eqnarray}
By Theorem~\ref{thm:2.1}, we have 
\begin{eqnarray}\label{eq9:thm:5.2}
\mathfrak{P}_{\bm{w},q}(\mathcal{A}^p_1,\cdots,\mathcal{A}^p_k) \preceq \lambda_{\min}^{p-1}K_1 K_2^{1/q}
\mathfrak{P}_{\bm{w},q}(\mathcal{A}_1,\cdots,\mathcal{A}_k),
\end{eqnarray}
then, we apply Lemma 3 in~\cite{chang2023randomI} to Eq.~\eqref{eq9:thm:5.2} to obtain Eq.~\eqref{eq1:thm:5.2}.
$\hfill \Box$

Following Theorem is about the deformed mean $\mathfrak{M}_{\sigma_{1/p}}$ and this theorem can be treated as the reverse counterparts of Lemma~\ref{lma:thm:4.1}.
\begin{theorem}\label{thm:5.3}
Let $\mathcal{A}_i \in \mathbb{C}^{I_1 \times \dots \times I_N \times I_1 \times \dots \times I_N}$ be random PD tensors satisfying $m\mathcal{I} \preceq \mathcal{A}_i \preceq M\mathcal{I}$ almost surely, a deterministic PD tensor $\mathcal{C} \in \mathbb{C}^{I_1 \times \dots \times I_N \times I_1 \times \dots \times I_N}$ and $\mathcal{X} \define \mathfrak{M}_{\sigma}(\mathcal{A}_1,\cdots,\mathcal{A}_k)$. Then, for any $p \geq 1$ and $r \geq 1$, we have
\begin{eqnarray}\label{eq1:thm:5.3}
\lefteqn{\mathrm{Pr}\left(\mathfrak{M}_{\sigma_{1/p}}(\mathcal{A}^p_1,\cdots,\mathcal{A}^p_k) \npreceq \mathcal{C}\right) \leq} \nonumber \\
&& \mathrm{Tr}\left(\mathbb{E}\left[\left(K_1\lambda^{p-1}_{\min}(\mathfrak{M}_{\sigma}(\mathcal{A}_1,\cdots,\mathcal{A}_k)) \mathfrak{M}_{\sigma}(\mathcal{A}_1,\cdots,\mathcal{A}_k)\right)^r\right] \star_N \mathcal{C}^{-1}\right),
\end{eqnarray}
where $K_1 \define K(M,m,p)$. On the other hand, we also have
\begin{eqnarray}\label{eq2:thm:5.3}
\lefteqn{\mathrm{Pr}\left(K_1^{-1}\lambda^{p-1}_{\max}(\mathfrak{M}_{\sigma}(\mathcal{A}_1,\cdots,\mathcal{A}_k)) \mathfrak{M}_{\sigma}(\mathcal{A}_1,\cdots,\mathcal{A}_k)\npreceq \mathcal{C}\right)
\leq} \nonumber \\
&& \mathrm{Tr}\left(\mathbb{E}\left[\left(
\mathfrak{M}_{\sigma_{1/p}}(\mathcal{A}^p_1,\cdots,\mathcal{A}^p_k)\right)^r\right] \star_N \mathcal{C}^{-1}\right).
\end{eqnarray}
\end{theorem}
\textbf{Proof:}
Because of $\lambda_{\min}(\mathcal{X}) \mathcal{X}^{-1} \preceq \mathcal{I}$ and Lemma~\ref{lma:5.1}, we have 
\begin{eqnarray}\label{eq3:thm:5.3}
\lefteqn{(\lambda^{1/2}_{\min}(\mathcal{X})\mathcal{X}^{-1/2}) \star_N \mathcal{A}_i^p \star_N (\lambda^{1/2}_{\min}(\mathcal{X})\mathcal{X}^{-1/2}) \preceq} \nonumber \\
&&K(M,m,p)\left[(\lambda^{1/2}_{\min}(\mathcal{X})\mathcal{X}^{-1/2}) \star_N \mathcal{A}_i \star_N (\lambda^{1/2}_{\min}(\mathcal{X})\mathcal{X}^{-1/2})\right]^p.
\end{eqnarray}
Similar to the proof in Lemma~\ref{lma:thm:4.1}, we will obtain
\begin{eqnarray}\label{eq4:thm:5.3}
\lefteqn{K^{-1}(M,m,p)\lambda^{p-1}_{\max}(\mathfrak{M}_{\sigma}(\mathcal{A}_1,\cdots,\mathcal{A}_k))\mathfrak{M}_{\sigma}(\mathcal{A}_1,\cdots,\mathcal{A}_k)} \nonumber \\
&\preceq& \mathfrak{M}_{\sigma_{1/p}}(\mathcal{A}^p_1,\cdots,\mathcal{A}^p_k) \nonumber \\
&\preceq& K(M,m,p)\lambda^{p-1}_{\min}(\mathfrak{M}_{\sigma}(\mathcal{A}_1,\cdots,\mathcal{A}_k))\mathfrak{M}_{\sigma}(\mathcal{A}_1,\cdots,\mathcal{A}_k).
\end{eqnarray}
Then, we can apply Lemma 3 in~\cite{chang2023randomI} to both inequalities Eq.~\eqref{eq4:thm:5.3} to obtain results of this theorem.
$\hfill \Box$

\bibliographystyle{IEEETran}
\bibliography{AHType_MultiArgu_Bib}

\begin{thebibliography}{10}
\providecommand{\url}[1]{#1}
\csname url@samestyle\endcsname
\providecommand{\newblock}{\relax}
\providecommand{\bibinfo}[2]{#2}
\providecommand{\BIBentrySTDinterwordspacing}{\spaceskip=0pt\relax}
\providecommand{\BIBentryALTinterwordstretchfactor}{4}
\providecommand{\BIBentryALTinterwordspacing}{\spaceskip=\fontdimen2\font plus
\BIBentryALTinterwordstretchfactor\fontdimen3\font minus
  \fontdimen4\font\relax}
\providecommand{\BIBforeignlanguage}[2]{{%
\expandafter\ifx\csname l@#1\endcsname\relax
\typeout{** WARNING: IEEEtran.bst: No hyphenation pattern has been}%
\typeout{** loaded for the language `#1'. Using the pattern for}%
\typeout{** the default language instead.}%
\else
\language=\csname l@#1\endcsname
\fi
#2}}
\providecommand{\BIBdecl}{\relax}
\BIBdecl

\bibitem{chang2023personalized}
S.~Y. Chang, H.-C. Wu, K.~Yan, X.~Chen, S.~C.-H. Huang, and Y.~Wu,
  ``Personalized multimedia recommendation systems using higher-order tensor
  singular-value-decomposition,'' \emph{IEEE Transactions on Broadcasting},
  2023.

\bibitem{chang2023TLS}
S.-Y. Chang and H.-C. Wu, ``Tensor-based least-squares solutions for
  multirelational signals and applications,'' \emph{IEEE Transactions on
  Cybernetics}, 2023.

\bibitem{chang2022tPI}
S.~Y. Chang and Y.~Wei, ``T-square tensors—part {I}: inequalities,''
  \emph{Computational and Applied Mathematics}, vol.~41, no.~1, p.~62, 2022.

\bibitem{chang2022tPII}
------, ``T-product tensors—part {II}: tail bounds for sums of random
  t-product tensors,'' \emph{Computational and Applied Mathematics}, vol.~41,
  no.~3, p.~99, 2022.

\bibitem{chang2022generaltail}
------, ``General tail bounds for random tensors summation: majorization
  approach,'' \emph{Journal of Computational and Applied Mathematics}, vol.
  416, p. 114533, 2022.

\bibitem{chang2022TKF}
S.~Y. Chang and H.-C. Wu, ``Tensor kalman filter and its applications,''
  \emph{IEEE Transactions on Knowledge and Data Engineering}, 2022.

\bibitem{chang2022tensorq}
------, ``Tensor quantization: High-dimensional data compression,'' \emph{IEEE
  Transactions on Circuits and Systems for Video Technology}, vol.~32, no.~8,
  pp. 5566--5580, 2022.

\bibitem{chang2022randouble}
S.~Y. Chang, ``Random double tensors integrals,'' \emph{arXiv preprint
  arXiv:2204.01927}, 2022.

\bibitem{chang2022TWF}
S.~Y. Chang and H.-C. Wu, ``Tensor wiener filter,'' \emph{IEEE Transactions on
  Signal Processing}, vol.~70, pp. 410--422, 2022.

\bibitem{chang2021convenient}
S.~Y. Chang and W.-W. Lin, ``Convenient tail bounds for sums of random
  tensors,'' \emph{Taiwanese Journal of Mathematics}, vol.~1, no.~1, pp. 1--36,
  2021.

\bibitem{chang2021tensorExp}
S.~Y. Chang, ``Tensor expander chernoff bounds,'' \emph{arXiv preprint
  arXiv:2105.06471}, 2021.

\bibitem{chang2022generalized}
------, ``Generalized hanson-wright inequality for random tensors,''
  \emph{arXiv preprint arXiv:2203.00659}, 2022.

\bibitem{chang2023randomI}
S.-Y. Chang, ``Random tensor inequalities and tail bounds for bivariate random
  tensor means, part i,'' \emph{arXiv preprint arXiv:2305.03301}, 2023.

\bibitem{chang2023randomII}
------, ``Random tensor inequalities and tail bounds for bivariate random
  tensor means, part ii,'' \emph{arXiv preprint arXiv:2305.03305}, 2023.

\bibitem{ando2004geometric}
T.~Ando, C.-K. Li, and R.~Mathias, ``Geometric means,'' \emph{Linear algebra
  and its applications}, vol. 385, pp. 305--334, 2004.

\bibitem{bhatia2006riemannian}
R.~Bhatia and J.~Holbrook, ``Riemannian geometry and matrix geometric means,''
  \emph{Linear algebra and its applications}, vol. 413, no. 2-3, pp. 594--618,
  2006.

\bibitem{palfia2016operator}
M.~P{\'a}lfia, ``Operator means of probability measures and generalized karcher
  equations,'' \emph{Advances in Mathematics}, vol. 289, pp. 951--1007, 2016.

\bibitem{lawson2014karcher}
J.~Lawson and Y.~Lim, ``Karcher means and karcher equations of positive
  definite operators,'' \emph{Transactions of the American Mathematical
  Society, Series B}, vol.~1, no.~1, pp. 1--22, 2014.

\bibitem{hiai2019ando}
F.~Hiai, Y.~Seo, and S.~Wada, ``Ando--hiai type inequalities for multivariate
  operator means,'' \emph{Linear and Multilinear Algebra}, vol.~67, no.~11, pp.
  2253--2281, 2019.

\bibitem{hiai2020ando}
------, ``{A}ndo--{H}iai-type inequalities for operator means and operator
  perspectives,'' \emph{International Journal of Mathematics}, vol.~31, no.~01,
  p. 2050007, 2020.

\bibitem{ni2019hermitian}
G.~Ni, ``Hermitian tensor and quantum mixed state,'' \emph{arXiv preprint
  arXiv:1902.02640}, 2019.

\bibitem{thompson1963certain}
A.~C. Thompson, ``On certain contraction mappings in a partially ordered vector
  space.'' \emph{Proceedings of the American Mathematical Society}, vol.~14,
  no.~3, pp. 438--443, 1963.

\bibitem{kubo1980means}
F.~Kubo and T.~Ando, ``Means of positive linear operators,''
  \emph{Mathematische Annalen}, vol. 246, pp. 205--224, 1980.

\bibitem{pales1997inverse}
Z.~P{\'a}les, ``Inverse and implicit function theorems for nonsmooth maps in
  banach spaces,'' \emph{Journal of Mathematical Analysis and Applications},
  vol. 209, no.~1, pp. 202--220, 1997.

\bibitem{ekeland1974variational}
I.~Ekeland, ``On the variational principle,'' \emph{Journal of Mathematical
  Analysis and Applications}, vol.~47, no.~2, pp. 324--353, 1974.

\bibitem{josip2013MondPecaricI}
J.~Pecaric, T.~Furuta, J.~M. Hot, and Y.~Seo, \emph{Mond-Pecaric method in
  operator inequalities: Inequalities for bounded selfadjoint operators on a
  Hilbert space, I}, 2005.

\end{thebibliography}

\end{document}